# Endogenous Aggregation of Multiple Data Envelopment Analysis Scores for Large Data Sets

**Hashem Omrani**[*]

Schulich School of Business, York University, Toronto, Ontario, Canada,

E-mail: omrani@yorku.ca

**Raha Imanirad**

Schulich School of Business, York University, Toronto, Ontario, Canada,

E-mail: raharad@schulich.yorku.ca

**Adam Diamant**

Schulich School of Business, York University, Toronto, Ontario, Canada,

E-mail: adiamant@schulich.yorku.ca

**Utkarsh Verma**

Schulich School of Business, York University, Toronto, Ontario, Canada,

E-mail: utkarshv@schulich.yorku.ca

**Amol Verma**

Department of Medicine, University of Toronto, Toronto, Ontario, Canada,

E-mail: amol.verma@mail.utoronto.ca

**Fahad Razak**

Department of Medicine, University of Toronto, Toronto, Ontario, Canada,

E-mail: fahad.razak@mail.utoronto.ca

[*] Corresponding Author, Email address: omrani@yorku.ca; omrani57@gmail.com

**Abstract**

We propose an approach for dynamic efficiency evaluation across multiple organizational dimensions using data envelopment analysis (DEA). The method generates both dimension-specific and aggregate efficiency scores, incorporates desirable and undesirable outputs, and is suitable for large-scale problem settings. Two regularized DEA models are introduced: a slack-based measure (SBM) and a linearized version of a nonlinear goal programming model (GP-SBM). While SBM estimates an aggregate efficiency score and then derives dimension-level efficiencies from the optimal slacks, GP-SBM first estimates dimension-level efficiencies and then derives an aggregate score. Both models utilize a regularization parameter to enhance discriminatory power while also directly integrating both desirable and undesirable outputs. We demonstrate the computational efficiency and validity of our approach on multiple datasets and apply it to a case study of twelve hospitals in Ontario, Canada, evaluating three theoretically grounded dimensions of organizational effectiveness over a 24-month period from January 2018 to December 2019: technical efficiency, clinical efficiency, and patient experience. Our numerical results show that SBM and GP-SBM better capture correlations among input/output variables and outperform conventional benchmarking methods that separately evaluate dimensions before aggregation.

## 1. Introduction

Data envelopment analysis (DEA) is a non-parametric method for estimating the efficiency of decision-making units (DMUs) with multiple input and output measures (Charnes et al., 1978; Mergoni et al., 2025). In recent years, advanced approaches have emerged, such as dynamic DEA models for evaluating DMUs over time (Fukuyama et al., 2023), and variable selection methods for applying DEA to large-scale systems (Lee and Cai, 2020). However, relatively little progress has been made in developing models that can simultaneously measure and integrate multiple, potentially correlated efficiency dimensions while also producing an aggregate score that is consistent with those dimensions. Indeed, most prior work evaluates each measure of a DMU's efficiency in isolation (Mahdiloo et al., 2015; Roth et al., 2019; Ferreira et al., 2023). Yet efficiency dimensions can conflict and overlapping input and output variables may introduce correlations that are not explicitly accounted for when the scores are jointly evaluated by decision makers.

When assessing organizations across multiple dimensions, some outputs are desirable while others may be deemed undesirable. Desirable outputs are those to be maximized, whereas undesirable outputs should be minimized to improve efficiency (Hua and Bian, 2007; Cheng and Zervopoulos, 2014). Incorporating both types of output variables into a unified model is a challenging but essential task for rigorous and comprehensive performance evaluation. Moreover, evaluating multiple dimensions of efficiency typically requires many input and output variables, which can reduce the discriminatory power of DEA-based models. To address this, methods such as Principal Component Analysis (PCA) (Adler and Golany, 2001) and the efficiency contribution measure (Pastor et al., 2002) have been proposed. However, these approaches are not well-suited to dynamic, multi-dimensional contexts, where a given dimension may depend on a single output that cannot be excluded or combined with other factors.

To address these challenges, we introduce two dynamic DEA formulations: a slack-based model (SBM) and a goal programming SBM (GP-SBM). SBM adopts a top-down approach, first optimizing overall efficiency and then distributing the aggregate value across dimensions. In contrast, GP-SBM follows a bottom-up approach: it simultaneously optimizes dimension-level efficiency through a goal programming model and then derives an aggregate efficiency score from the optimal solution. This distinction provides policymakers with complementary approaches that emphasize overall performance maximization (SBM) or fairness and balance across multiple organizational dimensions (GP-SBM). Both models incorporate undesirable outputs without data transformations (Cheng and Zervopoulos, 2014; Halkos and Petrou, 2019) and introduce a regularization term (Qin and Song, 2014) to enhance discriminatory power in settings with many input and output variables. While SBM is formulated as a linear program, GP-SBM is inherently nonlinear. Thus, we introduce a linear approximation of GP-SBM that provides an upper bound on the associated efficiency scores. We also prove that the regularized versions SBM and GP-SBM yield efficiency scores that are less than or equal to those of their non-regularized counterparts.

We evaluate the proposed models across three numerical studies. The first constructs a stylized example with three dimensions and two periods, where optimal efficiency scores can be computed directly. This example demonstrates that our dynamic DEA models produce consistent and theoretically intuitive results. The second study applies the models to a supply chain management problem of sustainable supplier selection in a static setting with three *independent* dimensions (Mahdiloo et al., 2015). Here, we highlight the computational advantages of SBM and the linearized GP-SBM, while confirming that the proposed models produce consistent results. The third study examines a unique healthcare setting using a unique dataset of twelve hospitals evaluated over a 24-month period across three dimensions of organizational effectiveness – technical efficiency, clinical efficiency, and patient experience (Roth et al., 2019) – with overlapping input and output variables across dimensions and the inclusion of carry-overs. As compared to two multiple-criteria decision-making methods used for benchmarking purposes, our framework generates DMU rankings that more effectively account for the shared input and output variables across the different dimensions, leading to results with distinct interpretations and policy implications. Consequently, although DEA models have been widely used to assess hospital efficiency (Kao et al., 2021; Afonso et al., 2023), our approach enables a more holistic evaluation of these complex organizations.

The paper is organized as follows. Section 2 reviews the literature and highlights the research gap. Section 3 presents the mathematical formulation of SBM and GP-SBM as well as their regularized versions.

Section 4 presents the numerical results for the stylized example and multi-dimensional supply chain management problem. Section 5 describes the data and findings from the hospital case study, while Section 6 concludes. Additional data tables and all proofs are provided in the supplementary material.

## 2. Literature Review and Positioning

In this section, we review prior research on dynamic DEA models in (i) integrating multiple efficiency dimensions; (ii) handling undesirable outputs; and (iii) enhancing discriminatory power in large datasets.

### 2.1. Integrating Multiple Dimensions of Efficiency

In practice, performance evaluation often involves several dimensions that must be assessed simultaneously. For example, banks measure both operational and profitability efficiency (Paradi et al., 2011), eco-efficiency studies integrate technical and environmental factors (Mahdiloo et al., 2015), and sustainability assessments encompass technical, social, and environmental dimensions (Omrani et al., 2022). A central challenge lies in integrating these diverse perspectives within a unified framework.

Previous studies have proposed two main approaches. The first requires input and output variables to be categorized in advance, with the decision maker specifying which contribute to each efficiency dimension instead of allowing them to be determined endogenously (Paradi et al., 2011). The second combines multiple dimensions into a single aggregate measure using goal programming within a multiplicative DEA framework (Mahdiloo et al., 2015; Omrani et al., 2022). However, this approach does not permit the calculation of separate efficiency scores for individual dimensions. When undesirable outputs are present, envelopment-based models such as SBM are preferable, as they directly incorporate slacks and effectively handle undesirable outputs (Cheng and Zervopoulos, 2014; Halkos and Petrou, 2019). Our proposed models build on these strengths by operating in envelopment form, accounting for slacks, and generating both aggregate and dimension-specific efficiency scores within a single optimization model.

In healthcare systems, hospitals seek to improve performance across multiple dimensions, yet most studies report a single aggregated score (Sommersguter-Reichmann, 2022). Afonso et al. (2023) link hospital efficiency with health quality and access to care using a network DEA model, while Ortega-Díaz and Martín (2022) examine the trade-off between hospital efficiency and service quality by estimating SBM models with and without quality measures as inputs. Although informative, these approaches estimate separate models for different efficiency dimensions, ignoring potential correlations among shared input and output variables. For instance, Roth et al. (2019) define technical efficiency, clinical efficiency, and patient experience in U.S. hospitals, but independently estimate each dimension. While this separation yields useful managerial insights, it overlooks potential interdependencies between the metrics. In contrast, our dynamic DEA models address this limitation by simultaneously evaluating multiple dimensions within a unified framework, capturing both complementarities and trade-offs. This integrated perspective offers a broader and more practical efficiency assessment than previous approaches. For comprehensive reviews of healthcare efficiency studies, see Cantor and Poh (2018) and Pai et al. (2024).

### 2.2. Handling Undesirable Outputs

One of the main challenges in DEA is determining how to incorporate undesirable outputs – such as $CO_2$ emissions, the number of defective products, or mortality rates – into the modeling framework. A common approach is to convert undesirable variables into desirable ones (Sommersguter-Reichmann, 2022). Notable examples include the linear additive transformation (Ali and Seiford, 1990) and the nonlinear multiplicative transformation (Lovell et al., 1995). These methods are widely used in practice (e.g., Lin et al., 2019; Roth et al., 2019), however, modifying the original data structure may lead to biased or unreliable results (Halkos and Petrou, 2019). Other approaches incorporate undesirable outputs by treating them as inputs to the DEA model (Liu and Sharp, 1999) or by imposing restrictive assumptions, such as requiring that undesirable outputs decrease only in the same proportion as the desirable outputs (Fare et al., 1989). One common technique uses directional distance functions (Chung et al., 1997), which adjust outputs along a fixed direction vector (Halkos and Petrou, 2019). However, restricting changes to a single direction parameter reduces model flexibility and complicates efficiency estimation across multiple dimensions. In contrast, slack-based models incorporate undesirable outputs directly and allow DMUs to address slacks by reducing inputs (input-oriented), increasing outputs (output-oriented), or adjusting both simultaneously (Tone, 2001). This flexibility makes slack-based models especially suitable for DMUs with mixed input and output variables, including undesirable ones (Cheng and Zervopoulos, 2014; Emrouznejad et al., 2025), and enables both dimension-specific and aggregate-level adjustments.

Slack-based models have been widely applied in healthcare to assess operational efficiency, including hospital performance evaluation (Zarrin, 2023), resource allocation and productivity monitoring (Gong et al., 2023; Sun et al., 2023), and health reform assessment (Zhang et al., 2018). These models allow input and output variables to be adjusted independently (Emrouznejad et al., 2025), which is particularly important in healthcare settings where undesirable performance indicators such as mortality rates, readmission rates, and patient waiting times often require distinct interventions for improvement, as their costs and resource intensity differ substantially. Most prior research has relied on static models that ignore temporal dynamics and evaluate only a single efficiency dimension. Moreover, existing models cannot endogenously estimate multiple efficiency dimensions while retaining strong discriminatory power for large datasets. Our proposed models address these limitations.

### 2.3. Applying DEA Models to Large-Scale Problem Settings

As the number of input and output variables increases, DEA models lose discriminatory power, with many DMUs receiving an efficiency score of one (Dyson et al., 2001; Lee and Cai, 2020). Several methods address this by proposing weight restrictions, dimensionality reduction, and super-efficiency extensions (e.g., Chen et al., 2016). Weight-based models – such as weight restrictions, assurance region models, cone ratio models, and common-weight models – limit DEA's flexibility in assigning weights (Omrani, 2013). However, as these methods use the multiplicative form of DEA, their applicability to envelopment-based models that incorporate undesirable outputs is limited. Another line of research emphasizes data reduction and variable selection. Common techniques include PCA (Adler and Golany, 2001; Peixoto et al., 2020), efficiency-contribution measures (Pastor et al., 2002), regression-based tests (Ruggiero, 2005), bootstrapping (Simar and Wilson, 2001), and average efficiency changes (Wagner and Shimshak, 2007). Empirical evidence suggests that efficiency-contribution and regression-based methods are effective for large datasets with highly correlated variables (Nataraja and Johnson, 2011; Chen and Han, 2021), whereas PCA with DEA suits smaller or more uncertain data structures (e.g., Ferreira et al., 2023). Lee and Cai (2020) develop a sign-constrained convex nonparametric least squares (SCNLS) model that uses LASSO for variable selection, while Benitez-Pena et al. (2020) propose a mixed-integer linear programming model that fixes the number of input and output variables to maximize efficiency. Wang et al. (2025, 2026) explore non-negative least squares and LASSO-type approaches for variable selection in high-dimensional DEA.

While many of the above methods can enhance the discriminatory power of DEA models, they remain confined to static, multiplicative DEA settings that estimate only a single efficiency measure. As such, they are not directly applicable to dynamic, multi-dimensional settings with undesirable outputs. Furthermore, these approaches generally assume that the primary objective is to reduce the total number of input and/or output variables, an assumption that is inappropriate in our context, where eliminating variables risks discarding essential information. Accordingly, we introduce a regularization parameter into our dynamic DEA models to enhance discriminatory power without removing key variables from any dimension. While the proposed methodology is broadly applicable to sectors such as banking, energy, and sustainability, we focus on healthcare, where large longitudinal and multidimensional datasets are common.

### 2.4. Research Gap

While an extensive body of literature proposes novel DEA models and evaluates hospital efficiency, several limitations remain. As shown in the structured overview provided in Supplementary Material 1, previous methodological contributions to multidimensional DEA typically focus on one of the following aspects: improving discriminatory power (Adler & Golany, 2001; Pastor et al., 2002), incorporating dynamic structures (Tone & Tsutsui, 2010), or evaluating multiple efficiency dimensions (Mahdiloo et al., 2015; Omrani et al., 2022). However, these aspects are typically addressed in isolation rather than integrated into a unified modeling framework. In particular, no previous study simultaneously considers multiple efficiency dimensions while accounting for the bias introduced by shared input/output variables. Another important issue relates to the discriminatory power of DEA models. Existing approaches – such as variable reduction or screening methods – primarily target settings where the total number of variables is large relative to the number of DMUs. In multidimensional settings, however, the number of variables within a specific dimension may be small, even if the aggregate count across all dimensions is high. Consequently, existing methodologies are not well suited to address the loss of discriminatory power inherent in such multidimensional structures.

Supplementary Material 1 also summarizes recent applications of DEA to hospital efficiency evaluation. As demonstrated, most studies employ traditional DEA, SBM, or directional distance function (DDF) models to estimate a single efficiency measure. Some studies analyze multiple aspects of hospital performance, but these dimensions are typically evaluated independently rather than jointly. Because dimension-specific efficiency scores are often correlated through shared variables, independent evaluations may lead to biased rankings. Our framework explicitly accounts for this structural interdependence. Furthermore, while a few studies seek to enhance discriminatory power in hospital efficiency analysis, they typically rely on methods such as super-efficiency, PCA, or multi-criteria decision-making, which do not address the endogenous weighting of dimensions.

Consequently, there is no DEA framework that simultaneously: (i) evaluates multiple efficiency dimensions; (ii) incorporates both desirable and undesirable input and output variables; (iii) enhances discriminatory power for dimensions with few variables; and (iv) endogenously derives both dimension-specific and aggregate efficiency measures. Our proposed framework addresses this gap by introducing two multidimensional DEA models and their regularized versions.

## 3. Methodology

This section presents the proposed models. First, we introduce the linear dynamic SBM and GP-SBM models, followed by their regularized versions. The models' characteristics are summarized below:

***Set Sizes and Index Definitions:***

- $I$: Number of inputs ($i$)
- $n$: Number of DMUs
- $R$: Number of desirable outputs ($r$)
- $K$: Number of undesirable outputs (k)
- $o$ Index of the DMU under evaluation
- $C$: Number of desirable carry-overs ($c$)
- $F$: Number of undesirable carry-overs ($f$)
- $G$: Number of efficiency dimensions ($g$)
- $T$: Number of periods ($t$)

***Data Parameters:***

$x_{ijt}$ : The value of the *i*-th input for the *j*-th DMU in period *t*

$y_{rjtg}$ : The value of the *r*-th desirable output for the *j*-th DMU in period *t,* used for measuring the *g-th* efficiency dimension

$b_{kjtg}$ : The value of the *k*-th undesirable output for the *j*-th DMU in period *t,* used for measuring the *g-th* efficiency dimension

$e_{cjt}$ : The value of the *c*-th desirable carry-over for the *j*-th DMU in period *t*

$z_{fjt}$ : The value of the *f*-th undesirable carry-over for the *j*-th DMU in period *t*

$\rho_{it}^{x}$ : The value of the regularized parameter for the *i*-th input in period *t*

$\rho_{rtg}^{y}, \rho_{ktg}^{b}$ : The value of the regularized parameter for the *r*-th desirable and the *k*-th undesirable outputs in period *t,* used for measuring the *g-th* efficiency dimension

$\rho_{ct}^{e}, \rho_{ft}^{z}$ : The value of the regularization parameter for *c*-th desirable and the *f*-th undesirable carry-overs in period *t,* used for measuring the *g-th* efficiency dimension

***Decision variables:***

$\lambda_{jt}, \Lambda_{jt}$ : Intensity variables associated with the *j*-th DMU in period *t*

$s_{it}^{x}, \overline{s}_{it}^{x}$ : Slack variables related to the *i*-th input in period *t*

$s_{rtg}^{y}, \overline{s}_{rtg}^{y}$ : Slack variable related to the *r*-th desirable output in period *t,* used for measuring the *g-th* efficiency dimension

$s_{ktg}^{b}, \overline{s}_{ktg}^{b}$ : Slack variables related to the *k*-th undesirable output in period *t,* used for measuring the *g-th* efficiency dimension

$s_{ct}^{e}, \overline{s}_{ct}^{e}$ : Slack variables related to the *c*-th desirable carry-over in period *t*

$s_{ft}^{z}, \overline{s}_{ft}^{z}$ : Slack variables related to the *f*-th undesirable carry-over in period *t*

$d_{gt}^{-}, d_{gt}^{+}$ : The negative and positive deviation variables of the *g-th* efficiency dimension

$v_i$ : The weight of *i*-th input variable

$u_{rg}, \mu_{kg}$ : Weight associated with the *r*-th desirable and the *k*-th undesirable outputs, used for measuring the *g-th* efficiency dimension

$w_c, \gamma_f$ : Weight associated with the *c*-th desirable and the *f*-th undesirable carry-overs, used for measuring the *g-th* efficiency dimension

We now propose two multidimensional DEA approaches that simultaneously evaluate dimension-specific and aggregate efficiency dimensions while accounting for desirable and undesirable indicators. An overview of the methodological workflow of the models is illustrated in Supplementary Material 2.

### 3.1. Dynamic SBM

Assume DMU *j* (*j*=1,..,*n*) consumes inputs $x_{ijt}$ (*i*=1,…,*I*; *t*=1,…,*T*) to produce desirable outputs $y_{rjtg}$ (*r*=1,…,*R*; *g*=1,…,*G*) while simultaneously generating undesirable outputs $b_{kjtg}$ (*k*=1,…,*K*) in period *t*. In addition, there are desirable carry-overs $e_{cjt}$ ($c$=1,…,$C$) and undesirable carry-overs $z_{fjt}$ ($f$=1,…,$F$) between periods $t$ and $t$+1. Suppose there are $G \geq 1$ efficiency dimensions. Following Tone and Tsutsui (2010), we formulate a variable-returns-to-scale (VRS) version of the dynamic SBM model, which we extend to simultaneously evaluate $G$ distinct efficiency dimensions plus an overall aggregate measure of efficiency.

$$
\begin{aligned}
&\min_{\lambda_{jt}, s_{it}^{x}, s_{rtg}^{y}, s_{ktg}^{b}, s_{ct}^{e}, s_{ft}^{z} \geq 0} \frac{\frac{1}{T}\sum_{t=1}^{T}\left[1-\frac{1}{I+F}[\sum_{i=1}^{I}\frac{s_{it}^{x}}{x_{iot}}+\sum_{f=1}^{F}\frac{s_{ft}^{z}}{z_{fot}}]\right]}{\frac{1}{T}\sum_{t=1}^{T}\left[1+\frac{1}{R+K+C}\sum_{g=1}^{G}[\sum_{i=1}^{I}\frac{s_{rtg}^{y}}{y_{rotg}}+\sum_{k=1}^{K}\frac{s_{ktg}^{b}}{b_{kotg}}+\sum_{c=1}^{C}\frac{s_{ct}^{e}}{e_{cot}}]\right]} \quad s.t: \\
&\sum_{j=1}^{n}\lambda_{jt}x_{ijt}+s_{it}^{x}=x_{iot}, \qquad i=1,...,I, t=1,...T \\
&\sum_{j=1}^{n}\lambda_{jt}y_{rjtg}-s_{rtg}^{y}=y_{rotg}, \qquad r=1,...,R, g=1,...,G, t=1,...,T \\
&\sum_{j=1}^{n}\lambda_{jt}b_{kjtg}+s_{ktg}^{b}=b_{kotg}, \qquad k=1,...,K, g=1,...,G, t=1,...,T \\
&\sum_{j=1}^{n}\lambda_{jt}e_{cjt}-s_{ct}^{e}=e_{cot}, \qquad c=1,...,C, t=1,...,T \\
&\sum_{j=1}^{n}\lambda_{jt}z_{fjt}+s_{ft}^{z}=z_{fot}, \qquad f=1,...,F, t=1,...,T \\
&\sum_{j=1}^{n}\lambda_{jt}e_{cjt}=\sum_{j=1}^{n}\lambda_{j(t+1)}e_{cjt}, \qquad c=1,...,C, t=1,...,T-1 \\
&\sum_{j=1}^{n}\lambda_{jt}z_{fjt}=\sum_{j=1}^{n}\lambda_{j(t+1)}z_{fjt}, \qquad f=1,...,F, t=1,...,T-1 \\
&\sum_{j=1}^{n}\lambda_{jt}=1, \qquad t=1,...,T
\end{aligned} \tag{1}
$$

The objective function maximizes overall efficiency by optimizing all slack variables. The first constraint reduces input consumption by maximizing input slacks. The second and third constraints correspond to desirable and undesirable outputs, respectively; the model seeks to increase desirable outputs while reducing undesirable ones. The fourth and fifth constraints handle carry-overs by maximizing desirable and minimizing undesirable components. Inter-period continuity of desirable and undesirable carry-overs between $t$ and $t$ + 1 is enforced by the sixth and seventh constraints. If the constraint $\sum_{j=1}^{n}\lambda_{jt}=1$ is removed from (1), the model reduces to a constant-returns-to-scale (CRS) form of the dynamic SBM. Note that if $G = 1$, (1) reduces to the dynamic SBM model in Tone and Tsutsui (2010). However, when $G > 1$, there are additional constraints and the denominator in the objective function is more complex. The following lemma reformulates the linear fractional optimization problem into an equivalent linear program.

***Lemma 1:*** Model (1) can be reformulated into a linear program as follows:

$$
\min_{q, \Lambda_{jt}, \overline{s}_{it}^{x}, \overline{s}_{rtg}^{y}, \overline{s}_{ktg}^{b}, \overline{s}_{ct}^{e}, \overline{s}_{ft}^{z} \geq 0} q-\frac{1}{T\times(I+F)}\left[\sum_{t=1}^{T}\sum_{i=1}^{I}\frac{\overline{s}_{it}^{x}}{x_{iot}}+\sum_{t=1}^{T}\sum_{f=1}^{F}\frac{\overline{s}_{ft}^{z}}{z_{fot}}\right] \quad s.t:
$$

$$
\begin{aligned}
& q+\frac{1}{T\times(R+K+C)}\sum_{t=1}^{T}\sum_{g=1}^{G}[\sum_{r=1}^{R}\frac{\bar{s}_{rtg}^{y}}{y_{rotg}}+\sum_{k=1}^{K}\frac{\bar{s}_{ktg}^{b}}{b_{kotg}}+\sum_{c=1}^{C}\frac{\bar{s}_{ct}^{e}}{e_{cot}}]=1 \\
& \sum_{j=1}^{n}\Lambda_{jt}x_{ijt}+\bar{s}_{it}^{x}=qx_{iot}, \qquad i=1,...,I, t=1,...T \\
& \sum_{j=1}^{n}\Lambda_{jt}y_{rjtg}-\bar{s}_{rtg}^{y}=qy_{rotg}, \qquad r=1,...,R, g=1,...,G, t=1,...,T \\
& \sum_{j=1}^{n}\Lambda_{jt}b_{kjtg}+\bar{s}_{ktg}^{b}=qb_{kotg}, \qquad k=1,...,K, g=1,...,G, t=1,...,T \\
& \sum_{j=1}^{n}\Lambda_{jt}e_{cjt}-\bar{s}_{ct}^{e}=qe_{cot}, \qquad c=1,...,C, t=1,...,T \\
& \sum_{j=1}^{n}\Lambda_{jt}z_{fjt}+\bar{s}_{ft}^{z}=qz_{fot}, \qquad f=1,...,F, t=1,...,T \\
& \sum_{j=1}^{n}\Lambda_{jt}e_{cjt}=\sum_{j=1}^{n}\Lambda_{j(t+1)}e_{cjt}, \qquad c=1,...,C, t=1,...,T-1 \\
& \sum_{j=1}^{n}\Lambda_{jt}z_{fjt}=\sum_{j=1}^{n}\Lambda_{j(t+1)}z_{fjt}, \qquad f=1,...,F, t=1,...,T-1 \\
& \sum_{j=1}^{n}\Lambda_{jt}=q, \qquad t=1,...,T
\end{aligned}
\tag{2}
$$

Let ( $q^{*},\Lambda_{jt}^{*},\bar{s}_{it}^{x*},\bar{s}_{rtg}^{y*},\bar{s}_{ktg}^{b*},\bar{s}_{ct}^{e*},\bar{s}_{ft}^{z*}$ ) be the optimal solution to the linear SBM model in (2). Then, the overall (aggregate) efficiency score of the *o*-th DMU under evaluation in period *t* is calculated as

$$
\begin{aligned}
\theta_{Overall,ot} &= q^{*}-\frac{1}{I+F}[\sum_{i=1}^{I}\frac{\bar{s}_{it}^{x*}}{x_{iot}}+\sum_{f=1}^{F}\frac{\bar{s}_{ft}^{z*}}{z_{fot}}]=\frac{q^{*}-\frac{1}{I+F}[\sum_{i=1}^{I}\frac{\bar{s}_{it}^{x*}}{x_{iot}}+\sum_{f=1}^{F}\frac{\bar{s}_{ft}^{z*}}{z_{fot}}]}{q^{*}+\frac{1}{R+K+C}\sum_{g=1}^{G}[\sum_{r=1}^{R}\frac{\bar{s}_{rtg}^{y*}}{y_{rotg}}+\sum_{k=1}^{K}\frac{\bar{s}_{ktg}^{b*}}{b_{kotg}}+\sum_{c=1}^{C}\frac{\bar{s}_{ct}^{e*}}{e_{cot}}]} \\
&= \frac{1-\frac{1}{I+F}[\sum_{i=1}^{I}\frac{s_{it}^{x*}}{x_{iot}}+\sum_{f=1}^{F}\frac{s_{ft}^{z*}}{z_{fot}}]}{1+\frac{1}{R+K+C}\sum_{g=1}^{G}[\sum_{r=1}^{R}\frac{s_{rtg}^{y*}}{y_{rotg}}+\sum_{k=1}^{K}\frac{s_{ktg}^{b*}}{b_{kotg}}+\sum_{c=1}^{C}\frac{s_{ct}^{e*}}{e_{cot}}]}
\end{aligned}
\tag{3}
$$

Thus, the mean overall efficiency of the *o*-th DMU over all time periods is given by $\theta_{Overall,o}=\frac{1}{T}\sum_{t=1}^{T}\theta_{Overall,ot}$ .

Using the optimal solution ( $q^{*},\Lambda_{jt}^{*},\bar{s}_{it}^{x*},\bar{s}_{rtg}^{y*},\bar{s}_{ktg}^{b*},\bar{s}_{ct}^{e*},\bar{s}_{ft}^{z*}$ ), we can also calculate the DEA score of the *g-th* efficiency dimension of the DMU under evaluation in period *t*, as[1]

---

[1] The directional distance function (DDF), is one of the most important approaches to handling undesirable outputs. However, in case studies involving multi-dimensional efficiency analysis, the DDF framework is not suitable if the goal is to obtain separate efficiency scores for each dimension. The DDF model produces equal scores for all dimensions when some specific directions have been selected. Consider the DDF model as follows: $\max\ \beta$ , s.t: $\sum_{j=1}^{n}\lambda_j x_{ij}\le x_{io}-\beta gx_i$ , $\sum_{j=1}^{n}\lambda_j y_{rj}\ge y_{ro}+\beta gy_r$ , $\sum_{j=1}^{n}\lambda_j b_{kj}=b_{ko}-\beta gb_k$ . Choosing some important directions such as

$$\theta_{got} = \frac{q^* - \frac{1}{I}\sum_{i=1}^{I}\frac{\bar{s}_{it}^{x*}}{x_{iot}}}{q^* + \frac{1}{R+K}[\sum_{r=1}^{R}\frac{\bar{s}_{rtg}^{y*}}{y_{rotg}} + \sum_{k=1}^{K}\frac{\bar{s}_{ktg}^{b*}}{b_{kotg}}]} = \frac{1 - \frac{1}{I}\sum_{i=1}^{I}\frac{s_{it}^{x*}}{x_{iot}}}{1 + \frac{1}{R+K}[\sum_{r=1}^{R}\frac{s_{rtg}^{y*}}{y_{rotg}} + \sum_{k=1}^{K}\frac{s_{ktg}^{b*}}{b_{kotg}}]} \tag{4}$$

Equation (4) builds on Tone and Tsutsui (2010), whose dynamic SBM framework yields only a single overall efficiency measure. In contrast, our formulation generalizes their approach by defining dimension-specific efficiencies derived from the optimal slacks of the aggregate solution.

In sum, the proposed model can be regarded as a top-down approach. That is, overall efficiency is first maximized after which the *g*-th dimension of efficiency is calculated using the resulting optimal solution. This means the model first constructs an efficiency frontier using all inputs, outputs, and carry-over variables collectively, thereby ensuring that each DMU's overall performance is assessed in an integrated manner. Once aggregate efficiency has been determined, the optimal slacks are applied to separately compute each efficiency dimension. Thus, the model derives dimension-specific efficiencies as outcomes of the overall optimal solution, consistent with a top-down evaluation approach. As indicated in (3) and (4), carry-overs contribute only to the measurement of overall efficiency. Finally, when all input, output, and carry-over slacks are zero, both the overall and dimension-specific efficiencies equal 1.0, which is the maximum attainable value. Conversely, an efficiency score of 1.0 implies the absence of any slacks.

### 3.2. Dynamic GP-SBM

Rather than defining a single objective function, we introduce multiple objectives (one for each efficiency dimension) and seek to optimize them together. More specifically, consider

$$\begin{aligned} &\min\left\{\frac{\frac{1}{T}\sum_{t=1}^{T}\left[1-\frac{1}{I}\sum_{i=1}^{I}\frac{s_{it}^{x}}{x_{iot}}\right]}{\frac{1}{T}\sum_{t=1}^{T}\left[1+\frac{1}{R+K}[\sum_{r=1}^{R}\frac{s_{rt,1}^{y}}{y_{rot,1}}+\sum_{k=1}^{K}\frac{s_{kt,1}^{b}}{b_{kot,1}}]\right]},\ldots,\frac{\frac{1}{T}\sum_{t=1}^{T}\left[1-\frac{1}{I}\sum_{i=1}^{I}\frac{s_{it}^{x}}{x_{iot}}\right]}{\frac{1}{T}\sum_{t=1}^{T}\left[1+\frac{1}{R+K}[\sum_{r=1}^{R}\frac{s_{rtG}^{y}}{y_{rotG}}+\sum_{k=1}^{K}\frac{s_{ktG}^{b}}{b_{kotG}}]\right]}\right\} \\ &= \min_{g=1,\ldots,G}\left\{\frac{1-\frac{1}{I\times T}\sum_{t=1}^{T}\sum_{i=1}^{I}\frac{s_{it}^{x}}{x_{iot}}}{1+\frac{1}{T\times(R+K)}\left[\sum_{t=1}^{T}\sum_{r=1}^{R}\frac{s_{rtg}^{y}}{y_{rotg}}+\sum_{t=1}^{T}\sum_{k=1}^{K}\frac{s_{ktg}^{b}}{b_{kotg}}\right]}\right\} \end{aligned} \tag{5}$$

---

$g = (gx_i, gy_r, gb_k) = (x_{iot}, y_{rotg}, b_{kotg})$ results in equal scores for all efficiency dimensions. In this case, all efficiency dimensions can be expressed as $\theta_{got} = \frac{1-\frac{1}{I}\sum_{i=1}^{I}\frac{\beta^* gx_i}{x_{iot}}}{1+\frac{1}{R+K}[\sum_{r=1}^{R}\frac{\beta^* gy_r}{y_{rotg}}+\sum_{k=1}^{K}\frac{\beta^* gb_k}{b_{kotg}}]} = \frac{1-\frac{1}{I}\sum_{i=1}^{I}\beta^*}{1+\frac{1}{R+K}[\sum_{r=1}^{R}\beta^*+\sum_{k=1}^{K}\beta^*]} = \frac{1-\beta^*}{1+\beta^*}$.

By including the same constraint set as in (1), (5) represents a multi-objective dynamic SBM formulation in which $G$ dimension-specific efficiencies are simultaneously minimized. To transform the multiple-objective, decision-making (MODM) model into a single-objective form, various techniques can be applied, including the parametric method (Saaty and Gass, 1954), the ε-constraint method (Mavrotas, 2009), and goal programming (Charnes et al., 1955). Since all $G$ objectives are nonlinear fractional functions, any direct approach – such as the weighted-sum method – remains nonlinear because each objective corresponds to a separate fractional measure, and the combined multi-objective formulation cannot be expressed as a single ratio. Consequently, the Charnes–Cooper transformation cannot be applied in this context to obtain a simplified model (Charnes and Cooper, 1962).

To address this limitation, we construct a linear approximation that provides an upper bound on the true efficiency score, allowing reformulation as a linear goal programming (GP) model. In DEA, full efficiency is defined as a score of 1.0, which naturally serves as the target in a GP framework. Thus, the GP version of the $g$-th objective in (5) can be written as $\dfrac{1-\dfrac{1}{I\times T}\sum_{t=1}^{T}\sum_{i=1}^{I}\dfrac{s_{it}^{x}}{x_{iot}}}{1+\dfrac{1}{T\times(R+K)}\left[\sum_{t=1}^{T}\sum_{r=1}^{R}\dfrac{s_{rtg}^{y}}{y_{rotg}}+\sum_{t=1}^{T}\sum_{k=1}^{K}\dfrac{s_{ktg}^{b}}{b_{kotg}}\right]}+d_g=1$,

where $d_g$ denotes the deviation variable. If $d_g(t)=0$, then all slacks are zero, and the $g$-th efficiency dimension equals 1. If $d_g(t)>0$, then at least one slack is positive, implying an efficiency below 1.0. Directly incorporating this goal constraint into the optimization problem introduces nonlinearity because the deviation terms interact with the ratio structure. To address this, we present the following result.

***Proposition 1:*** The efficiency scores in (5) are bounded from above by a linear function such that

$$\theta_{go}=\frac{1-\dfrac{1}{I\times T}\sum_{t=1}^{T}\sum_{i=1}^{I}\dfrac{s_{it}^{x}}{x_{iot}}}{1+\dfrac{1}{T\times(R+K)}\left[\sum_{t=1}^{T}\sum_{r=1}^{R}\dfrac{s_{rtg}^{y}}{y_{rotg}}+\sum_{t=1}^{T}\sum_{k=1}^{K}\dfrac{s_{ktg}^{b}}{b_{kotg}}\right]}\leq\frac{1}{I\times T}\sum_{t=1}^{T}\sum_{i=1}^{I}\frac{s_{it}^{x}}{x_{iot}}+\frac{1}{T\times(R+K)}\left[\sum_{t=1}^{T}\sum_{r=1}^{R}\frac{s_{rtg}^{y}}{y_{rotg}}+\sum_{t=1}^{T}\sum_{k=1}^{K}\frac{s_{ktg}^{b}}{b_{kotg}}\right]:=\tilde{\theta}_{go}\qquad(6)$$

Based on the relationship derived in (6), we construct a GP model with the linear goal constraint $\tilde{\theta}_{go}=d_{go}$. Although this represents an approximation of the true efficiency score, its accuracy relative to $\theta_{go}$ is evaluated numerically in Section 4 where we observe strong correspondence. In general, there are three standard techniques for solving GP problems: Lexicographic GP (LGP), Weighted GP (WGP), and Min–Max (Chebyshev) GP (Chang, 2007). The LGP approach requires the decision maker to pre-specify the priority levels of goals while WGP and Chebyshev GP do not. Accordingly, due to its simplicity, the weighted goal programming SBM model is formulated as follows:

$$
\begin{aligned}
&\min_{\lambda_{jt}, s_{it}^{x}, s_{rtg}^{y}, s_{ktg}^{b}, s_{ct}^{e}, s_{ft}^{z}, d^{-}, d^{+} \geq 0} \sum_{g=1}^{G} \left( -d_g \right), \quad s.t: \\
&\sum_{j=1}^{n} \lambda_{jt} x_{ijt} + s_{it}^{x} = x_{iot}, \qquad i = 1,...,I, t = 1,...T \\
&\sum_{j=1}^{n} \lambda_{jt} y_{rjtg} - s_{rtg}^{y} = y_{rotg}, \qquad r = 1,...,R, g = 1,...,G, t = 1,...,T \\
&\sum_{j=1}^{n} \lambda_{jt} b_{kjtg} + s_{ktg}^{b} = b_{kotg}, \qquad k = 1,...,K, g = 1,...,G, t = 1,...,T \\
&\sum_{j=1}^{n} \lambda_{jt} e_{cjt} - s_{ct}^{e} = e_{cot}, \qquad c = 1,...,C, t = 1,...,T \\
&\sum_{j=1}^{n} \lambda_{jt} z_{fjt} + s_{ft}^{z} = z_{fot}, \qquad f = 1,...,F, t = 1,...,T \\
&\sum_{j=1}^{n} \lambda_{jt} e_{cjt} = \sum_{j=1}^{n} \lambda_{j(t+1)} e_{cjt}, \qquad c = 1,...,C, t = 1,...,T-1 \\
&\sum_{j=1}^{n} \lambda_{jt} z_{fjt} = \sum_{j=1}^{n} \lambda_{j(t+1)} z_{fjt}, \qquad f = 1,...,F, t = 1,...,T-1 \\
&\sum_{j=1}^{n} \lambda_{jt} = 1, \qquad t = 1,...,T \\
&\frac{1}{I \times T} \sum_{t=1}^{T} \sum_{i=1}^{I} \frac{s_{it}^{x}}{x_{iot}} + \frac{1}{T \times (R+K)} \left[ \sum_{t=1}^{T} \sum_{r=1}^{R} \frac{s_{rtg}^{y}}{y_{rotg}} + \sum_{t=1}^{T} \sum_{k=1}^{K} \frac{s_{ktg}^{b}}{b_{kotg}} \right] = d_g, \quad g = 1,...,G
\end{aligned}
\tag{7}
$$

Model (7) follows a bottom-up approach, optimizing all efficiency dimensions simultaneously rather than sequentially. In this model, each dimension's efficiency is treated as a separate optimization goal within the model. Using goal programming, the model balances trade-offs across these dimensions within a unified framework. After finding the optimal solution, the overall efficiency of each DMU is calculated as the product of its individual dimension scores. That is, let ($\lambda_{jt}^{*}, s_{it}^{x*}, s_{rtg}^{y*}, s_{ktg}^{b*}, s_{ct}^{e*}, s_{ft}^{z*}$) denote the optimal solution of the GP-SBM model in (7). The overall (aggregate) efficiency of the *o*-th DMU in period *t* is computed using (3), while the *g*-th dimension efficiency is derived from (4).

### 3.3. Regularized Dynamic SBM and GP-SBM

Several methods, such as PCA-DEA, ECM, and AEC, have been developed to enhance the discriminatory power of DEA models. One key assumption in these methods is that there are many input/output variables, and the objective is to reduce or eliminate them to increase discriminatory power. However, in multi-dimensional, dynamic DEA models – where separate efficiency scores are calculated for different performance dimensions – this assumption does not hold. Although the total number of variables across all dimensions may be large, each individual dimension often relies on only a small subset. In such cases, reducing or combining variables through PCA or similar approaches may eliminate or distort important variables for certain dimensions, making it impossible to accurately compute those efficiency scores.

To improve the discriminatory power of DEA-based models, the weights of input, output, and carry-over variables in the objective function can be adjusted by introducing a penalty term (Qin and Song, 2014). This promotes sparsity by encouraging smaller – or, in some cases, zero – weights for certain variables, thereby minimizing their influence on efficiency scores. Since the envelopment form of DEA does not explicitly involve weights, this idea is adapted so that its influence is applied indirectly through the dual formulation. Specifically, by penalizing the magnitude of variable weights in the multiplier form of DEA, we indirectly affect the right-hand side of the envelopment model (as expressed in its dual), thereby enhancing the model's ability to discriminate between DMUs without eliminating any variables.

To facilitate our presentation, assume that $v_i, u_{rg}, \mu_{kg}, w_c$, and $\gamma_f$ are, respectively, the weights of inputs, desirable outputs, undesirable outputs, desirable carry-overs, and undesirable carry-overs. To incorporate regularization into the SBM model, the dual form of (2) is required. For ease of exposition, we present the dual for a single period (*T=1*), and thus, the sixth and seventh constraints of (2), which link flows between consecutive periods, are omitted. Thus, the multiplier form of (2) for a single period is

$$
\begin{aligned}
&\max_{v_i, u_{rg}, \mu_{kg}, w_c, \gamma_f, \Omega, \Psi} \Psi, \quad s.t: \\
&-\sum_{i=1}^{I} v_i x_{ij} + \sum_{g=1}^{G}\sum_{r=1}^{R} u_{rg} y_{rjg} - \sum_{g=1}^{G}\sum_{k=1}^{K} \mu_{kg} b_{kjg} + \sum_{c=1}^{C} w_c e_{cj} - \sum_{f=1}^{F} \gamma_f z_{kj} + \Omega \le 0, \quad j = 1,...,n \\
&\sum_{i=1}^{I} v_i x_{io} - \sum_{g=1}^{G}\sum_{r=1}^{R} u_{rg} y_{rog} + \sum_{g=1}^{G}\sum_{k=1}^{K} \mu_{kg} b_{kog} - \sum_{c=1}^{C} w_c e_{co} + \sum_{f=1}^{F} \gamma_f z_{ko} - \Omega + \Psi \le 1 \\
&-v_i \le -\frac{1}{I+F}\frac{1}{x_{io}}, i = 1,...,I \\
&-u_{rg} + \frac{1}{R+K+C}\frac{\Psi}{y_{rog}} \le 0, r = 1,...,R, g = 1,...,G \\
&-\mu_{kg} + \frac{1}{R+K+C}\frac{\Psi}{b_{kog}} \le 0, k = 1,...,K, g = 1,...,G \\
&-w_c + \frac{1}{R+K+C}\frac{\Psi}{e_{co}} \le 0, c = 1,...,C \\
&-\gamma_f \le -\frac{1}{I+F}\frac{1}{z_{fo}}, f = 1,...,F
\end{aligned}
\tag{8}
$$

To add regularization, the objective function in (8) can be rewritten as:

$$
\max \Psi - \sum_{i=1}^{I} v_i \rho_i^x - \sum_{g=1}^{G}\sum_{r=1}^{R} u_{rg} \rho_{rg}^y - \sum_{g=1}^{G}\sum_{k=1}^{K} \mu_{kg} \rho_{kg}^b - \sum_{c=1}^{C} w_c \rho_c^e - \sum_{f=1}^{F} \gamma_f \rho_f^z \tag{9}
$$

where $\rho_i^x, \rho_r^y, \rho_k^b, \rho_c^e$ and $\rho_f^z$ are, respectively, the regularization parameters for inputs, desirable outputs, undesirable outputs, desirable carry-overs, and undesirable carry-overs. The regularized objective function (9) encourages the *o*-th DMU to give lower weight (or zero) to less important inputs, outputs, and carry-overs, effectively filtering out variables with minimal influence. Building on these insights, by formulating the dynamic DEA version of (8) for *T>1* and replacing its objective with the regularized form in (9), we obtain a regularized variant of the dynamic SBM model in (2), which can be expressed as:

$$
\begin{aligned}
&\min_{q,\Lambda_{jt},\overline{s}_{it}^{x},\overline{s}_{rtg}^{y},\overline{s}_{ktg}^{b},\overline{s}_{ct}^{e},\overline{s}_{ft}^{z}\geq 0} q-\frac{1}{T\times(I+F)}\left[\sum_{t=1}^{T}\sum_{i=1}^{I}\frac{\overline{s}_{it}^{x}}{x_{iot}}+\sum_{t=1}^{T}\sum_{f=1}^{F}\frac{\overline{s}_{ft}^{z}}{z_{fot}}\right], \quad s.t: \\
&q+\frac{1}{T\times(R+K+C)}\sum_{t=1}^{T}\sum_{g=1}^{G}[\sum_{r=1}^{R}\frac{\overline{s}_{rtg}^{y}}{y_{rotg}}+\sum_{k=1}^{K}\frac{\overline{s}_{ktg}^{b}}{b_{kotg}}+\sum_{c=1}^{C}\frac{\overline{s}_{ct}^{e}}{e_{cot}}]=1 \\
&\sum_{j=1}^{n}\Lambda_{jt}x_{ijt}+\overline{s}_{it}^{x}=qx_{iot}+\rho_{it}^{x}, \qquad i=1,\ldots,I,t=1,\ldots T \\
&\sum_{j=1}^{n}\Lambda_{jt}y_{rjtg}-\overline{s}_{rtg}^{y}=qy_{rotg}-\rho_{rtg}^{y}, \qquad r=1,\ldots,R,g=1,\ldots,G,t=1,\ldots,T \\
&\sum_{j=1}^{n}\Lambda_{jt}b_{kjtg}+\overline{s}_{ktg}^{b}=qb_{kotg}+\rho_{ktg}^{b}, \qquad k=1,\ldots,K,g=1,\ldots,G,t=1,\ldots,T \\
&\sum_{j=1}^{n}\Lambda_{jt}e_{cjt}-\overline{s}_{ct}^{e}=qe_{cot}-\rho_{ct}^{e}, \qquad c=1,\ldots,C,t=1,\ldots,T \\
&\sum_{j=1}^{n}\Lambda_{jt}z_{fjt}+\overline{s}_{ft}^{z}=qz_{fot}+\rho_{ft}^{z}, \qquad f=1,\ldots,F,t=1,\ldots,T \\
&\sum_{j=1}^{n}\Lambda_{jt}e_{cjt}=\sum_{j=1}^{n}\Lambda_{j(t+1)}e_{cjt}, \qquad c=1,\ldots,C,t=1,\ldots,T-1 \\
&\sum_{j=1}^{n}\Lambda_{jt}z_{fjt}=\sum_{j=1}^{n}\Lambda_{j(t+1)}z_{fjt}, \qquad f=1,\ldots,F,t=1,\ldots,T-1 \\
&\sum_{j=1}^{n}\Lambda_{jt}=q, \qquad t=1,\ldots,T
\end{aligned}
\tag{10}
$$

Deriving this result follows directly from linear programming duality theory in that adding penalty terms to the objective in the primal formulation modifies the right-hand side of the constraints in the dual. Consequently, introducing regularization adjusts the feasible region of the envelopment form without removing or altering any input, output, or carry-over variables. Further, in (10), the regularization parameters appear as positive terms on the right-hand side for inputs, undesirable outputs, and undesirable carry-overs, and as negative terms for desirable outputs and desirable carry-overs. This structure aligns with the intended direction of improvement: it penalizes excessive inputs, carry-overs, and undesirable outputs while discouraging reductions in desirable ones. Similarly, the regularized counterpart of (7) can be obtained by modifying the objective of the multiplier form to include regularization and obtaining its dual.

$$
\begin{aligned}
&\min_{\lambda_{jt},s_{it}^{x},s_{rtg}^{y},s_{ktg}^{b},s_{ct}^{e},s_{ft}^{z}\geq 0} \sum_{g=1}^{G}\left(-d_{g}\right), \quad s.t: \\
&\sum_{j=1}^{n}\lambda_{jt}x_{ijt}+s_{it}^{x}=x_{iot}+\rho_{it}^{x}, \qquad i=1,\ldots,I,t=1,\ldots T \\
&\sum_{j=1}^{n}\lambda_{jt}y_{rjtg}-s_{rtg}^{y}=y_{rotg}-\rho_{rtg}^{y}, \qquad r=1,\ldots,R,g=1,\ldots,G,t=1,\ldots,T
\end{aligned}
$$

$$\sum_{j=1}^{n} \lambda_{jt} b_{kjtg} + s_{ktg}^{b} = b_{kotg} + \rho_{ktg}^{b}, \qquad k = 1,...,K, g = 1,...,G, t = 1,...,T \tag{11}$$

$$\sum_{j=1}^{n} \lambda_{jt} e_{cjt} - s_{ct}^{e} = e_{cot} - \rho_{ct}^{e}, \qquad c = 1,...,C, t = 1,...,T$$

$$\sum_{j=1}^{n} \lambda_{jt} z_{fjt} + s_{ft}^{z} = z_{fot} + \rho_{ft}^{z}, \qquad f = 1,...,F, t = 1,...,T$$

$$\sum_{j=1}^{n} \lambda_{jt} e_{cjt} = \sum_{j=1}^{n} \lambda_{j(t+1)} e_{cjt}, \qquad c = 1,...,C, t = 1,...,T-1$$

$$\sum_{j=1}^{n} \lambda_{jt} z_{fjt} = \sum_{j=1}^{n} \lambda_{j(t+1)} z_{fjt}, \qquad f = 1,...,F, t = 1,...,T-1$$

$$\sum_{j=1}^{n} \lambda_{jt} = 1, \qquad t = 1,...,T$$

$$\frac{1}{I \times T} \sum_{t=1}^{T} \sum_{i=1}^{I} \frac{s_{it}^{x}}{x_{iot}} + \frac{1}{T \times (R+K)} \left[ \sum_{t=1}^{T} \sum_{r=1}^{R} \frac{s_{rtg}^{y}}{y_{rotg}} + \sum_{t=1}^{T} \sum_{k=1}^{K} \frac{s_{ktg}^{b}}{b_{kotg}} \right] = d_{g}, \quad g = 1,...,G$$

As previously noted, after models (10) and (11) are solved, the dimension-specific and overall efficiency scores can be obtained by applying (3) and (4). Furthermore, we note two interesting properties.

***Proposition 2*****:** The efficiency scores generated by the regularized models of SBM and GP-SBM are less than or equal to those produced by the non-regularized models.

***Proposition 3*****:** Let $\rho$ be a common regularization parameter for all variables. A sufficient condition to ensure the non-negativity of the efficiency scores is to choose

$$\rho \leq \min_{r,c,o,t,g} \left\{ q^{*} y_{rotg}, q^{*} e_{\text{cot}}, \frac{I - \sum_{i=1}^{I} \frac{s_{it}^{x*}}{x_{iot}}}{\sum_{i=1}^{I} \frac{\rho}{x_{iot}}}, \frac{I + F - [\sum_{i=1}^{I} \frac{s_{it}^{x*}}{x_{iot}} + \sum_{f=1}^{F} \frac{s_{ft}^{z*}}{z_{fot}}]}{[\sum_{i=1}^{I} \frac{1}{x_{iot}} + \sum_{f=1}^{F} \frac{1}{z_{fot}}]} \right\} \tag{12}$$

Since $\rho$ parametrizes the right-hand-side value of constraints (10) and (11), any user-defined value satisfying (12) ensures a non-empty feasible region, thereby yielding admissible efficiency scores.

***Proposition 4 (Monotonicity in ρ)***: If $0 \leq \rho_1 \leq \rho_2$, then $\theta(\rho_1) \geq \theta(\rho_2)$ where $\theta(.)$ denotes the efficiency score of a DMU.

## 4. Numerical Examples

This section presents two numerical examples illustrating the performance of the proposed dynamic SBM and GP-SBM models. The first is a simulated case designed to demonstrate the models' behavior in a dynamic, multi-dimensional setting. Because the true DMU rankings can be derived analytically, this example illustrates that the model's outputs align with theoretical expectations. The second example, adapted from the literature (Mahdiloo et al. 2015), highlights the computational advantages of the SBM and linearized GP-SBM models in a static setting (*T = 1*) with three independent dimensions.

Numerical Example 1: This example uses a small dataset designed to highlight the main features of our models. It includes four DMUs – A, B, C, and D – that use one input $x$ to produce one desirable output $y$, two undesirable outputs $b_1$ and $b_2$, and one undesirable carry-over $z$. The data are shown in Table 1. We evaluate three efficiency dimensions: one based on a desirable output and two based on undesirable outputs.

**Table 1: Data for the synthetic numerical example**

| | $T_1$ | | | | | $T_2$ | | | | |
|---|---|---|---|---|---|---|---|---|---|---|
| | $x$ | $y$ | $b_1$ | $b_2$ | $z$ | $x$ | $y$ | $b_1$ | $b_2$ | $z$ |
| ***A*** | 1 | 2 | 1 | 1 | 1 | 2 | 4 | 1 | 2 | 2 |
| ***B*** | 1 | 3 | 1 | 1 | 1 | 2 | 2 | 3 | 3 | 2 |
| ***C*** | 1 | 1 | 2 | 1 | 1 | 2 | 2 | 2 | 3 | 2 |
| ***D*** | 1 | 1 | 2 | 2 | 1 | 2 | 3 | 1 | 2 | 2 |

Let $\theta_{g,t}$ be the *g-th* efficiency dimension in period $t$ where: (i) $\theta_{1,t}$ denotes input $x$, desirable output y, and carry-over $z$; (ii) $\theta_{2,t}$ denotes input $x$, undesirable output $b_1$, and carry-over $z$; and (iii) $\theta_{3,t}$ denotes input $x$, undesirable output $b_2$, and carry-over $z$. This example verifies whether the models yield valid results when certain dimensions consist solely of undesirable factors. For example, in the second dimension, DMUs use the same input $x$ and carry-over $z$, but differ in the amount of undesirable output $b_1$ they produce. In this case, DMUs A and B generate lower levels of $b_1$ than C and D and are therefore expected to achieve higher efficiency scores. In other words, $A \simeq B \gg C \simeq D$, which means that the efficiency scores of DMUs A and B, and also C and D, are approximately equal. However, the efficiency scores of DMUs A and B should be greater than those of C and D. As a consequence, the rankings of the DMUs in time periods $T_1$ and $T_2$ can be derived as follows:

| | $T_1$ | $T_2$ |
|---|---|---|
| $\theta_{1,t}$: | $B \gg A \gg C \simeq D$ | $A \gg D \gg B \simeq C$ |
| $\theta_{2,t}$: | $A \simeq B \gg C \simeq D$ | $A \simeq D \gg C \gg B$ |
| $\theta_{3,t}$: | $A \simeq B \simeq C \gg D$ | $A \simeq D \gg B \simeq C$ |
| $\theta_{Overall,t}$: | $B \gg A \gg C \gg D$ | $A \gg D \gg C \gg B$ |

Table 2 reports the results for the regularized and non-regularized versions of the SBM model (10) and the GP-SBM model (11) for this example. The data was first normalized by dividing each value by the maximum value of its corresponding column within each period. For both models, the same regularization parameter was applied to all variables ($\rho = 0.03$). Values in parentheses indicate the efficiency scores from the non-regularized models. The results show that both the SBM and GP-SBM models accurately estimate overall and dimension-specific efficiencies, with the resulting scores aligning with the expected ranking patterns. In addition, the GP-SBM models generally produce slightly higher scores than the SBM models. These findings confirm that the proposed models produce consistent results in multi-dimensional settings, and in particular, a pathological case where certain dimensions contain only undesirable outputs.

**Table 2: Results of regularized SBM and GP-SBM models for the synthetic numerical example**

| Period | Efficiency | Models | A | B | C | D |
|---|---|---|---|---|---|---|
| ***T*₁** | $\theta_1$ | *SBM* | 0.622(0.667) | 0.919 (1.000) | 0.299(0.333) | 0.302(0.333) |
| | | *GP-SBM* | 0.628(0.667) | 0.942(1.000) | 0.314(0.333) | 0.314(0.333) |
| | $\theta_2$ | *SBM* | 0.903(1.000) | 0.884(1.000) | 0.608(0.667) | 0.613(0.667) |
| | | *GP-SBM* | 0.915(1.000) | 0.915(1.000) | 0.634(0.667) | 0.634(0.667) |
| | $\theta_3$ | *SBM* | 0.903(1.000) | 0.884(1.000) | 0.853(1.000) | 0.613(0.667) |
| | | *GP-SBM* | 0.915(1.000) | 0.915(1.000) | 0.915(1.000) | 0.634(0.667) |
| | $\theta_{\text{Overall}}$ | *SBM* | 0.785(0.857) | 0.896(1.000) | 0.487(0.545) | 0.457(0.50) |
| | | *GP-SBM* | 0.794(0.857) | 0.924(1.000) | 0.512(0.545) | 0.473(0.50) |
| ***T*₂** | $\theta_1$ | *SBM* | 0.933(1.000) | 0.460(0.500) | 0.449(0.500) | 0.679(0.750) |
| | | *GP-SBM* | 0.942(1.000) | 0.471(0.500) | 0.471(0.500) | 0.706(0.750) |
| | $\theta_2$ | *SBM* | 0.875(1.000) | 0.561(0.600) | 0.598(0.667) | 0.828(1.000) |
| | | *GP-SBM* | 0.890(1.000) | 0.572(0.600) | 0.628(0.667) | 0.890(1.000) |
| | $\theta_3$ | *SBM* | 0.918(1.000) | 0.697(0.750) | 0.681(0.750) | 0.885(1.000) |
| | | *GP-SBM* | 0.928(1.000) | 0.711(0.750) | 0.711(0.750) | 0.928(1.000) |
| | $\theta_{\text{Overall}}$ | *SBM* | 0.908(1.000) | 0.557(0.600) | 0.559(0.621) | 0.788(0.900) |
| | | *GP-SBM* | 0.919(1.000) | 0.568(0.600) | 0.586(0.621) | 0.829(0.900) |

*Numerical Example 2:* Because few studies have applied dynamic DEA models across multiple efficiency dimensions, we compare the proposed SBM and GP-SBM models using the dataset from Mahdiloo et al. (2015), who evaluate two independent constructs in sustainable supplier selection. Specifically, Mahdiloo et al. (2015) assess and rank 20 suppliers across two efficiency dimensions (technical and environmental) and derive an overall score termed "eco-efficiency". Their analysis employs a goal programming-based multiplier version of DEA with a single time-period to estimate both the individual dimension efficiencies and the combined eco-efficiency score. The input-output structure of their study is as follows. For technical efficiency, the inputs are the number of employees and energy consumption, while the outputs include sales, return on assets (ROA), and environmental R&D investment. For environmental efficiency, no explicit inputs are defined. Instead, this dimension evaluates how effectively suppliers generate desirable outputs (sales, ROA, and environmental R&D investment) relative to the undesirable output ($CO_2$ emissions). Finally, the eco-efficiency score integrates all variables to represent the suppliers' overall performance, rather than serving as an aggregate measure derived from the scores of each individual dimension.

To adapt the SBM and GP-SBM models to this static (single-period) setting, we remove the carry-over constraints. Due to the simplified formulation, we can compare the results of the linear (7) and nonlinear (5) versions of the GP-SBM model. To facilitate computational analysis, we apply the Charnes–Cooper transformation to the goal constraint $\dfrac{1-\dfrac{1}{I\times T}\sum_{t=1}^{T}\sum_{i=1}^{I}\dfrac{s_{it}^{x}}{x_{iot}}}{1+\dfrac{1}{T\times(R+K)}\left[\sum_{t=1}^{T}\sum_{r=1}^{R}\dfrac{s_{rtg}^{y}}{y_{rotg}}+\sum_{t=1}^{T}\sum_{k=1}^{K}\dfrac{s_{ktg}^{b}}{b_{kotg}}\right]}+d_g=1$. Denoting $q_g = \dfrac{1}{1+\dfrac{1}{T\times(R+K)}\left[\sum_{t=1}^{T}\sum_{r=1}^{R}\dfrac{s_{rtg}^{y}}{y_{rotg}}+\sum_{t=1}^{T}\sum_{k=1}^{K}\dfrac{s_{ktg}^{b}}{b_{kotg}}\right]}$, the goal constraint can be re-expressed as

$q_g + \dfrac{1}{T\times(R+K)}\left[\sum_{t=1}^{T}\sum_{r=1}^{R}\dfrac{q_g s_{rtg}^{y}}{y_{rotg}}+\sum_{t=1}^{T}\sum_{k=1}^{K}\dfrac{q_g s_{ktg}^{b}}{b_{kotg}}\right]=1$ which is bilinear. Thus, we replace the final constraint in (11) with the bilinear constraint above and refer to the resulting nonlinear GP-SBM model as NLP-GP.

The efficiency scores generated by the different models are presented in Supplementary Material 4. Table 3 reports the Spearman correlation coefficients between the ranks obtained from applying each of these models to the data. All correlations are highly significant at the 0.01 level. Notably, GP-SBM and NLP-GP show stronger consistency with Mahdiloo et al. (2015) (denoted as MO) than the SBM model across all three sustainability dimensions. Although the efficiency scores produced by GP-SBM and NLP-GP are remarkably similar, the NLP-GP model is bilinear and non-convex, meaning that convergence to a global optimum is not guaranteed. Furthermore, in our experiments, NLP-GP required substantially longer computation times to return a solution, rendering it impractical for large-scale problem instances.

**Table 3: Spearman correlation between models across different dimensions**

| | Technical | | | | Environmental | | | | Eco-Efficiency | | | |
|---|---|---|---|---|---|---|---|---|---|---|---|---|
| | SBM | GP-SBM | NLP-GP | MO | SBM | GP-SBM | NLP-GP | MO | SBM | GP-SBM | NLP-GP | MO |
| SBM | 1 | 0.946 | 0.946 | 0.784 | 1 | 0.989 | 0.989 | 0.841 | 1 | 0.944 | 0.944 | 0.827 |
| GP-SBM | | 1 | 1.000 | 0.900 | | 1 | 1.000 | 0.860 | | 1 | 1.000 | 0.904 |
| NLP-GP | | | 1 | 0.900 | | | 1 | 0.860 | | | 1 | 0.904 |
| MO | | | | 1 | | | | 1 | | | | 1 |

*** All correlations are significant at the 0.01 level

## 5. Case Study: Integrating Multiple Dimensions in Hospital Efficiency Assessment

We next apply the proposed dynamic SBM and GP-SBM models to a large dataset to evaluate multiple dimensions of hospital efficiency. Specifically, we estimate monthly efficiency scores for 12 hospitals in Ontario, Canada, over a 24-month period from January 2018 to December 2019. The analysis focuses on three theoretically grounded dimensions of hospital effectiveness (Roth et al., 2019): technical efficiency, clinical efficiency, and patient experience. This application is particularly appropriate for our framework, as hospital operations encompass multiple interdependent processes that share inputs and produce both desirable and undesirable outcomes. Previous studies have used a wide range of input and output variables to assess hospital performance (Cantor and Poh, 2018). Supplementary Material 5 summarizes the commonly used variables identified in the literature. In addition, Fazria and Dhamayanti (2021) conduct a comprehensive review of hospital efficiency studies and find that the most frequently used inputs such as the number of beds, medical personnel, operational costs, non-medical staff, and medical technician personnel. The most common outputs are the number of inpatients, outpatient visits, surgeries, days of inpatient care, and emergency visits.

Our analysis uses data collected from GEMINI (Verma et al. 2025), an organization that compiles inpatient data from multiple hospitals across Ontario, Canada. Guided by prior research and data availability, we selected the following variables for assessing technical, clinical, and patient experience:

*Input Variables:*

- Number of physicians: The total count of physicians providing medical care in the hospital.
- Bed occupancy levels: The percentage of available hospital beds occupied by patients.

*Desirable Outputs for Measuring Technical Efficiency:*

- Number of emergency visits: The total visits to the emergency room (ER).
- Number of surgeries: The count of in-hospital surgical procedures that are performed (the data does not include planned surgeries, as those data were not available in the database).
- Number of discharges: The total number of patients released from the hospital after treatment (this includes patients discharged home with support, to nursing homes, or to supported living facilities).

*Undesirable Output for Measuring Clinical Efficiency:*

- Mortality rate: The proportion of patients who die after being admitted to the hospital.

*Undesirable Outputs for Measuring Patient Experience:*

- Triage time in the ER: The average time patients spend after registering for emergency or ambulatory care until their first interaction with an ER physician (hours).
- Wait times for radiography/echocardiograms: The average waiting period (hours) for patients to receive imaging services (i.e., the time difference between the performed and ordered dates/times).

*Undesirable Carry-Overs:*

- 30-day readmission rate: The percentage of patients readmitted to the hospital within 30 days of being discharged from any hospital within Ontario.
- Length of hospital stay: The average duration of a patient's stay in the hospital (days).

All input and output variables were aggregated monthly. Furthermore, mortality rate, 30-day readmission rate, and length of hospital stay were risk-adjusted to account for local demographics according to the Ontario Health Guidelines.[2] For technical efficiency, hospitals consume resources such as physicians and beds to generate desirable outputs. However, due to data limitations, all outputs for clinical and patient experience efficiencies are undesirable. To provide additional clarity on the structure of the efficiency evaluation problem considered in this study, a conceptual illustration of the hospital efficiency framework is provided in Supplementary Material 2. The figure illustrates how hospitals use shared inputs to generate outcomes across three efficiency dimensions, while undesirable inputs and dimension-specific carryover variables link performance across periods. This setup captures both intra- and inter-dimensional linkages, with the regularization parameter facilitating their discrimination. These endogenous relationships are not accounted for in the existing multidimensional DEA literature, which typically assumes independence among inputs, outputs, and carryovers. Table 4 provides detailed descriptive statistics, including the minimum, maximum, mean, and standard deviation for each variable across the 24-month study period. All data are first normalized in each period (month) using the following formula (Yang et al., 2018):

$$\tilde{x}_{ijt} = \frac{x_{ijt}}{\max\limits_j\{x_{jt}\}}, \tilde{y}_{rjt} = \frac{y_{rjt}}{\max\limits_j\{y_{rjt}\}}, \tilde{b}_{kjt} = \frac{b_{kjt}}{\max\limits_j\{b_{kjt}\}}, \tilde{e}_{cjt} = \frac{e_{cjt}}{\max\limits_j\{e_{cjt}\}}, \tilde{z}_{fjt} = \frac{z_{fjt}}{\max\limits_j\{z_{fjt}\}}, \forall i,r,k,c,f,t \quad (13)$$

[2] *https://www.hqontario.ca/Portals/0/documents/qi/gemqin-opr-background-and-indicator-details.pdf*

Given the size of the dataset, our analysis focuses on the regularized versions of SBM and *linear* GP-SBM. This choice is motivated by the relatively large number of inputs and outputs (there are 10) as compared to the number of DMUs per period. Under these conditions, the non-regularized SBM and GP-SBM models assign an efficiency score of one to all hospitals. Furthermore, the nonlinear version of GP-SBM is computationally intractable in this big data setting. For simplicity, we assume that all regularization parameters for input, outputs, and carry-overs are the same. In other words, $\rho = \rho_{it}^x = \rho_{rt}^y = \rho_{kt}^b = \rho_{ct}^e = \rho_{ft}^z$. However, we perform a sensitivity analysis by varying ρ to assess its impact on the resulting DEA scores. Finally, as a benchmark, we separately estimate technical, clinical, and patient experience efficiencies and then aggregate the corresponding scores using the Technique for Order Preference by Similarity to Ideal Solution (TOPSIS) and Complex Proportional Assessment (COPRAS) methods.

**Table 4: A summary of raw data**

| | | | **Min** | **Max** | **Mean** | **SD** |
|---|---|---|---|---|---|---|
| | **Inputs** | Number of physicians | 20 | 125 | 57.7 | 25.41 |
| | | Bed occupancy level (%) | 62.1 | 96.2 | 86.8 | 5.41 |
| **Outputs** | **Technical Efficiency** | Number of emergency visits | 1.0 | 1013 | 91.47 | 115.78 |
| | | Number of surgeries | 1.0 | 222 | 31.94 | 51.21 |
| | | Number of discharges | 252 | 1052 | 621.3 | 229.2 |
| | **Clinical Efficiency** | Mortality rate (%) | 4.6 | 9.4 | 7.6 | 1.03 |
| | **Patient Experience** | Triage time in the ER (hours) | 0.001 | 0.201 | 0.073 | 0.063 |
| | | Wait time for radiography (hours) | 0.68 | 16.98 | 4.27 | 3.00 |
| | **Carry-Overs** | 30-day readmission rate (%) | 13.12 | 16.85 | 14.15 | 0.722 |
| | | Length of hospital stay (days) | 7.96 | 10.06 | 9.23 | 0.412 |

### 5.1. Results of Regularized SBM and Linear GP-SBM

Since the normalized data are less than or equal to one, the regularization parameter is intentionally kept small ($\rho = 0.005$). Table 5 reports the mean technical, clinical, patient experience, and overall efficiency scores over time; more detailed results are provided in Supplementary Materials 6 and 7.

For technical efficiency, the SBM model (10) yields scores ranging from 0.5433 (H9) to 0.9892 (H3), indicating that H3 is the most technically efficient hospital, while H9 exhibits the weakest performance. For the linear GP-SBM model (11), the range of technical efficiency scores is 0.6355 (H9) to 0.9917 (H3), reflecting generally higher efficiency estimates across hospitals. For example, the technical efficiency of H1 increases from 0.6787 under SBM to 0.7781 under linear GP-SBM, illustrating the tendency of linear GP-SBM to generate higher efficiency estimates. For clinical efficiency, the SBM model indicates that H8 is the best performing hospital (efficiency equal to 0.9887) and H9 is the worst (0.9230). Under the linear GP-SBM model, the range of clinical efficiency scores is from 0.9253 (H10) to 0.9894 (H8), which is higher as compared with the SBM model, showing overall improvements in clinical performance. Regarding patient experience efficiency, the SBM model shows that H8 is again the top-ranked hospital with a score of 0.9849, while H12 is the least efficient, with a score of 0.6144. In the linear GP-SBM model, the patient experience efficiency scores also improve over SBM, ranging from 0.6662 (H12) to 0.9856 (H8). Finally, for overall efficiency, the SBM model reports scores ranging from 0.6068 (H2) to 0.9644 (H3). Again, linear GP-SBM has higher values, ranging from 0.7202 (H2) to 0.9657 (H3). These results indicate that while hospitals such as H3 and H8 perform better across dimensions, other hospitals (e.g., H2, H9, H12) need to use resources more efficiently to improve their levels of performance.

**Table 5: The mean of efficiency scores generated by regularized SBM and linear GP-SBM models over time**

| | **Regularized SBM** | | | | **Regularized Linear GP-SBM** | | | |
|---|---|---|---|---|---|---|---|---|
| | *Technical* | *Clinical* | *Patient* | *Overall* | *Technical* | *Clinical* | *Patient* | *Overall* |
| **H1** | 0.6787 | 0.9737 | 0.7732 | 0.7159 | 0.7781 | 0.9839 | 0.7338 | 0.7758 |
| **H2** | 0.5689 | 0.9698 | 0.6453 | 0.6068 | 0.6898 | 0.9812 | 0.7135 | 0.7202 |
| **H3** | 0.9892 | 0.9856 | 0.9317 | 0.9644 | 0.9917 | 0.9864 | 0.9279 | 0.9657 |
| **H4** | 0.8669 | 0.9797 | 0.9698 | 0.9165 | 0.8813 | 0.9805 | 0.9644 | 0.9234 |
| **H5** | 0.7503 | 0.9848 | 0.9282 | 0.8314 | 0.7837 | 0.9870 | 0.9362 | 0.8554 |
| **H6** | 0.8014 | 0.9861 | 0.9615 | 0.8695 | 0.8244 | 0.9875 | 0.9638 | 0.8850 |
| **H7** | 0.9132 | 0.9626 | 0.9012 | 0.9131 | 0.9254 | 0.9647 | 0.9014 | 0.9199 |
| **H8** | 0.8792 | 0.9887 | 0.9849 | 0.9276 | 0.8867 | 0.9894 | 0.9856 | 0.9323 |
| **H9** | 0.5433 | 0.9230 | 0.9652 | 0.6814 | 0.6355 | 0.9435 | 0.9708 | 0.7586 |
| **H10** | 0.7627 | 0.9319 | 0.8385 | 0.8019 | 0.8078 | 0.9253 | 0.8188 | 0.8276 |
| **H11** | 0.9330 | 0.9826 | 0.9762 | 0.9535 | 0.9369 | 0.9834 | 0.9751 | 0.9554 |
| **H12** | 0.9070 | 0.9802 | 0.6144 | 0.7926 | 0.9241 | 0.9842 | 0.6662 | 0.8272 |

Figures 1 shows the efficiencies generated by the SBM and linear GP-SBM models for all periods, respectively. In the technical dimension, both SBM and linear GP-SBM show that efficiency varies widely across hospitals. Some hospitals, such as H3 and H11, consistently perform well, while others, including H1, H2, and H9, exhibit lower average values and greater variability in scores. As a result, H3 and H11 can be considered as benchmarks for all hospitals in this dimension. In addition, the results reveal that some hospitals face operational challenges, which should be taken into account for future improvement. In contrast, performance in the clinical dimension differs significantly. As shown in Figure 1, all hospitals demonstrate consistently high clinical efficiency (except H10) as the mortality rate remained low. Finally, hospital performance in the patient experience dimension is generally stronger than in the technical dimension. As shown in Figure 1, hospitals H1, H2, and H10 exhibit greater variability in performance as compared to the other institutions. Moreover, H2, H10 (under linear GP-SBM), and H12 record lower efficiency scores, suggesting that these hospitals face challenges related to patient wait times.

By endogenously aggregating technical, clinical, and patient experience efficiencies, we obtain an overall measure of hospital performance. In doing so, we observe that several hospitals achieve strong and balanced performance. For example, both SBM and the linear GP-SBM model show that hospitals H3, H4, H8, and H11 maintain high and stable overall efficiencies across all periods. In contrast, hospitals H1, H2, H9, and H10 display greater variability, with efficiency scores fluctuating considerably over time. Overall, these results suggest that while clinical efficiency remains consistently strong, patient experience performance is more variable, and technical and overall efficiencies show greater heterogeneity.

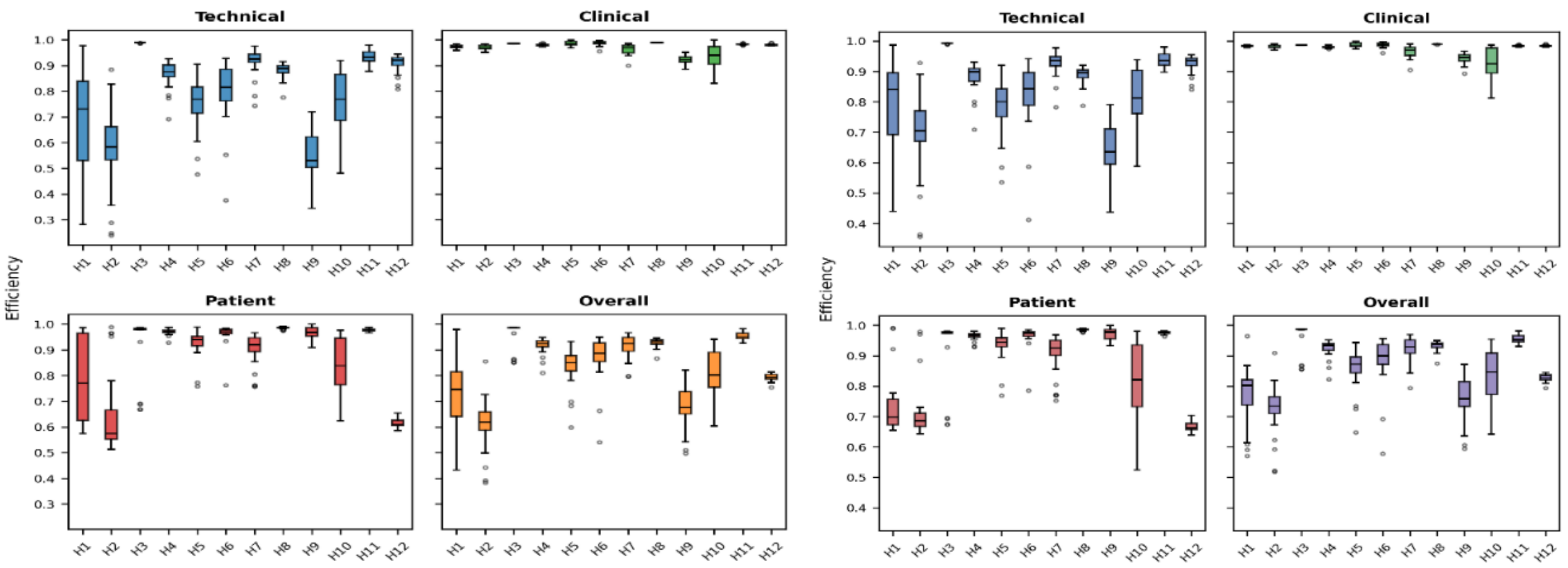

**Figure 1: Efficiency generated by SBM model (left) and GP-SBM (right) across multiple dimensions**

From a managerial perspective, the results suggest that hospital benchmarking should not rely on separately estimated efficiency dimensions. Because technical, clinical, and patient experience share common resources and processes, a hospital's relative ranking may change when these dimensions are jointly evaluated. This is particularly relevant for hospitals with unbalanced performance across dimensions, as a single input or output may disproportionately influence their efficiency computation.

### 5.2. Sensitivity Analysis and Robustness

In this section, we assess how variations in the value of the regularization parameter (ρ) affects the average technical, clinical, patient, and overall efficiencies estimated by the SBM and linear GP-SBM models. To this end, Figure 2 shows that efficiency gradually decrease as ρ increases, indicating improved discriminatory capability. Across all tested values of ρ, the linear GP-SBM model consistently achieves higher efficiency scores than SBM. Figure 2 further illustrates that the SBM model is more sensitive to variations in the regularization parameter (ρ). In the technical dimension, the SBM scores start at 1 when ρ = 0 and gradually decline to 0.644 at ρ = 0.01. In contrast, linear GP-SBM scores decrease from 1 (at ρ = 0) to 0.745 (at ρ = 0.01), demonstrating the model's robustness in maintaining higher efficiency under stronger regularization. When ρ = 0, both models yield an average efficiency of 1 across all dimensions, confirming that without regularization, they cannot effectively distinguish between hospitals.

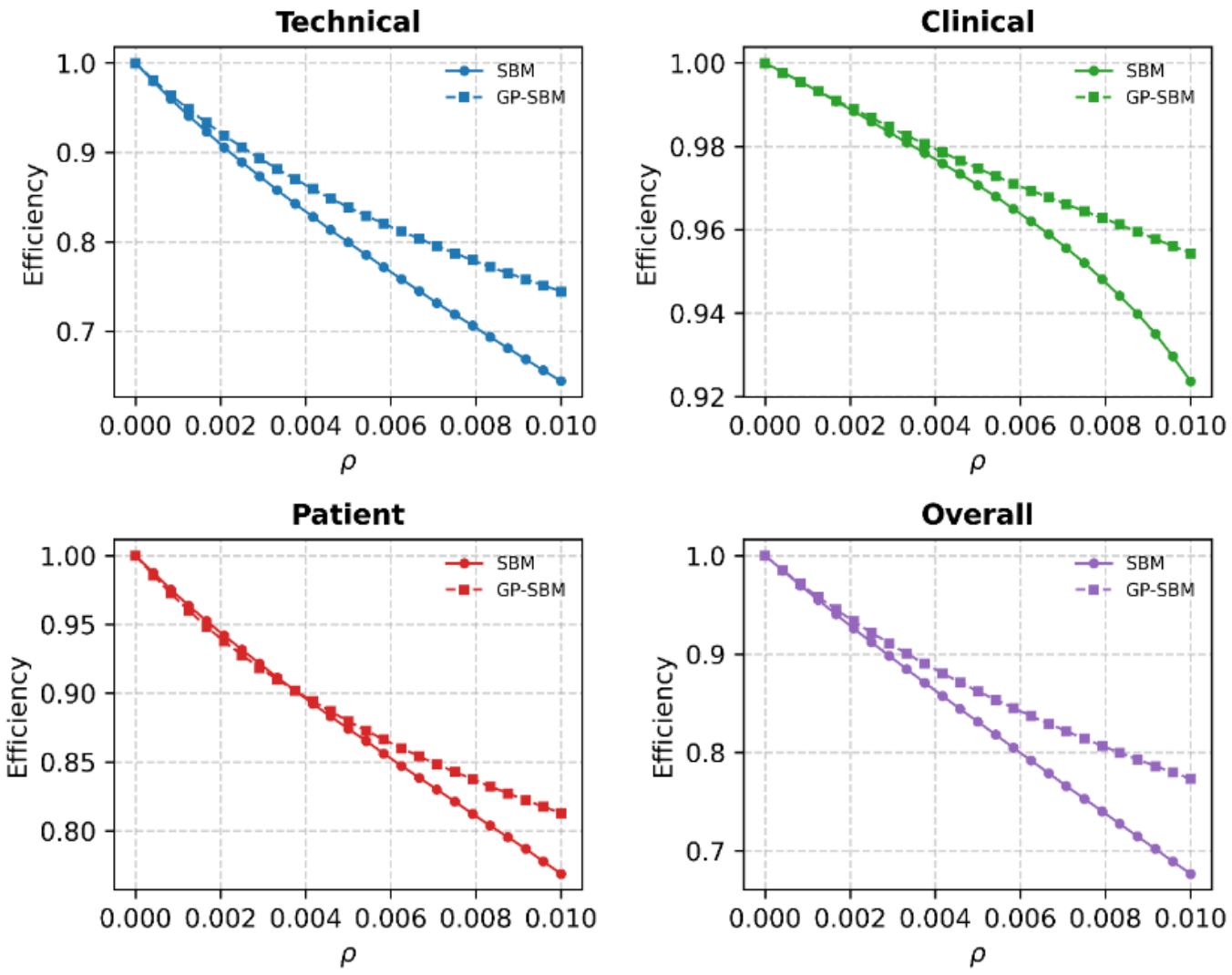

**Figure 2: Efficiency comparison of SBM and linear GP-SBM across the regularization parameter (ρ)**

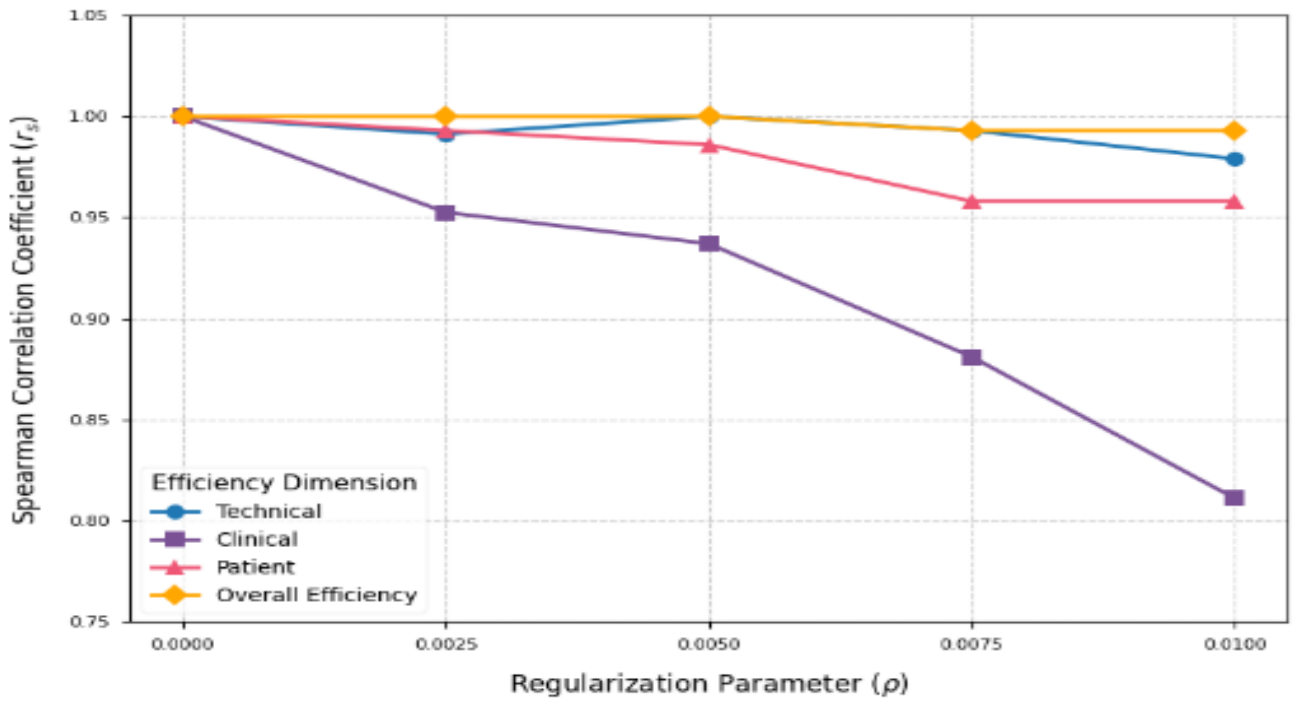

**Figure 3: Spearman correlation coefficients between SBM and GP-SBM for different ρ values**

Although the above analysis indicates that some differences are observed in the raw DEA efficiency scores, Table 6 indicates that the rankings produced by the SBM and linear GP-SBM models are highly consistent, particularly for the technical and overall efficiency dimensions. Specifically, Table 6 reports the Spearman correlation coefficients between the rankings generated by the two models for the average technical, clinical, patient experience, and overall efficiencies using the regularization parameter $\rho = 0.005$. The strong and statistically significant correlations indicate a high degree of concordance between the models across all efficiency dimensions. In the clinical and patient experience categories, all correlations are significant at the 0.01 level, underscoring the robustness of the relationship between the SBM and GP-SBM results. Moreover, the correlations across all dimensions and time periods remain statistically significant at least at the 0.1 level, confirming the consistency of model performance over time. Finally, Figure 3 confirms that the strong correlations among models persist even as the regularization parameter varies (see Supplementary Material 8 for a more detailed presentation of the correlation coefficients).

**Table 6: Spearman correlation between SBM and linear GP-SBM**

| ***Technical*** | SBM | GP-SBM | ***Clinical*** | SBM | GP-SBM |
|---|---|---|---|---|---|
| SBM | 1.000 | 1.000$^{***}$ | SBM | 1.000 | 0.937$^{***}$ |
| GP-SBM | | 1.000 | GP-SBM | | 1.000 |
| ***Patient experience*** | | | ***Overall*** | | |
| SBM | 1.000 | 0.986$^{***}$ | SBM | 1.000 | 1.000$^{***}$ |
| GP-SBM | | 1.000 | GP-SBM | | 1.000 |

$^{***}$ Correlation is significant at the 0.01 level

Intertemporal analysis shows strong persistence within efficiency dimensions, particularly for clinical and patient experience measures, where autocorrelations often exceed 0.90 ($p < 0.01$). In contrast, technical efficiency is more variable, indicating more frequent operational adjustments. Inter-dimensional effects are observed, however, as technical efficiency often serves as a leading indicator for patient outcomes (e.g., P4 → P5). These findings suggest that while hospital performance is primarily driven by intra-dimensional factors, dynamic inter-dimensional linkages do play a role. See Supplementary Material 9 for full results.

### 5.3. Validation and Benchmarking

To objectively validate the SBM and linear GP-SBM results, we compare the aggregate scores with those generated using multi-criteria decision-making (MCDM) methods, which are commonly used alongside DEA to integrate multiple efficiency dimensions into a single ranking (Omrani et al., 2022). Because the clinical dimension includes only one output and the patient experience dimension includes two, there is limited potential for dimension reduction (e.g., via PCA) or variable selection (e.g., efficiency-contribution or regression-based methods). Therefore, as a comparison, each efficiency dimension is evaluated independently and then aggregated using two established MCDM techniques: TOPSIS and COPRAS.

TOPSIS ranks DMUs based on their relative closeness to an ideal solution, whereas COPRAS evaluates alternatives by balancing beneficial and non-beneficial criteria. The integration of DEA with these methods is well established in the literature: TOPSIS has been used to combine efficiency dimensions in airline evaluation (Omrani et al., 2022) and to improve DEA's discriminatory power (Rakhshan, 2017), while COPRAS has been applied in domains such as financial portfolio selection (Gupta et al., 2019) and renewable energy assessment (Yilmaz, 2023). The procedural steps of TOPSIS are detailed in Behzadian et al. (2012), and those of COPRAS in Zavadskas et al. (2007). In both methods, equal weights are assigned to each efficiency dimension when computing the final score.

Table 7 summarizes the results. Columns 2-4 report average efficiency scores for each dimension, calculated independently of the others, while Columns 5-8 present the aggregate scores and corresponding ranks (in parentheses) obtained from TOPSIS, COPRAS, SBM, and linear GP-SBM. Unlike TOPSIS and COPRAS, which compute the final score by averaging the dimension-specific efficiencies from Columns 2-4, the SBM and linear GP-SBM models directly generate both individual and overall efficiency scores. As expected, hospitals H3 and H8 consistently achieve the highest rankings across all four methods, whereas H2 and H9 rank lowest. The primary differences occur among mid-ranked hospitals (e.g., H4, H5, H6, H7, H10, and H12), where rankings vary depending on the method applied. These discrepancies reflect fundamental methodological distinctions: SBM and linear GP-SBM, as DEA models, evaluate DMUs relative to an efficiency frontier, while TOPSIS and COPRAS are distance-based ranking techniques.

**Table 7: The results of the TOPSIS, COPRAS, SBM, and linear GP-SBM models**

| | Technical | Clinical | Patient | TOPSIS | COPRAS | SBM | GP-SBM |
|---|---|---|---|---|---|---|---|
| **H1** | 0.166 | 0.929 | 0.606 | 0.455 (9) | 0.555 (10) | 0.716 (10) | 0.776 (10) |
| **H2** | 0.106 | 0.922 | 0.607 | 0.438 (10) | 0.529 (11) | 0.607 (12) | 0.720 (12) |
| **H3** | 0.990 | 0.987 | 0.877 | 0.928 (2) | 1.000 (1) | 0.964 (1) | 0.966 (1) |
| **H4** | 0.530 | 0.739 | 0.835 | 0.553 (6) | 0.728 (7) | 0.916 (4) | 0.923 (4) |
| **H5** | 0.444 | 0.973 | 0.831 | 0.666 (4) | 0.760 (5) | 0.831 (7) | 0.855 (7) |
| **H6** | 0.796 | 0.974 | 0.908 | 0.869 (3) | 0.929 (3) | 0.869 (6) | 0.885 (6) |
| **H7** | 0.854 | 0.643 | 0.397 | 0.398 (11) | 0.673 (9) | 0.913 (5) | 0.919 (5) |
| **H8** | 0.878 | 0.979 | 0.949 | 0.928 (1) | 0.979 (2) | 0.928 (3) | 0.932 (3) |
| **H9** | 0.399 | 0.570 | 0.474 | 0.196 (12) | 0.497 (12) | 0.681 (11) | 0.759 (11) |
| **H10** | 0.570 | 0.771 | 0.652 | 0.490 (7) | 0.688 (8) | 0.802 (8) | 0.828 (8) |
| **H11** | 0.932 | 0.654 | 0.488 | 0.456 (8) | 0.740 (6) | 0.953 (2) | 0.955 (2) |
| **H12** | 0.669 | 0.960 | 0.597 | 0.617 (5) | 0.763 (4) | 0.793 (9) | 0.827 (9) |

Table 8 reports the Spearman correlations among TOPSIS, COPRAS, SBM, and linear GP-SBM. All correlations are statistically significant, with SBM and linear GP-SBM producing identical rankings. However, their rankings differ across certain periods, as shown in Supplementary Material 10. The strong correlations between COPRAS and both DEA-based models support the validity of our approach. However, notable differences do emerge between our models and the MCDM benchmarks. These discrepancies arise because TOPSIS and COPRAS do not account for shared input and output variables across efficiency dimensions, whereas the proposed SBM and GP-SBM models explicitly incorporate these interdependencies. Consequently, by construction, when input and/or output variables contribute to multiple efficiency dimensions, our models offer a more accurate assessment of efficiency rankings.

These differences are also relevant in practice. When rankings are obtained from separate evaluation and ex-post aggregation, hospitals may be compared against peer institutions that do not reflect the joint structure of shared inputs and outputs. In contrast, the proposed SBM and GP-SBM models identify benchmark positions within a unified frontier, which can lead to different conclusions regarding relative performance and improvement priorities. Therefore, the proposed framework provides decision makers with a more coherent basis for benchmarking hospitals when efficiency dimensions are interdependent.

**Table 8: Spearman correlation test between different models**

| | SBM | GP-SBM | TOPSIS | COPRAS |
|---|---|---|---|---|
| **SBM** | 1 | 1.000*** | 0.5734* | 0.7133*** |
| **GP-SBM** | | 1 | 0.5734* | 0.7133*** |
| **TOPSIS** | | | 1 | 0.9440*** |
| **COPRAS** | | | | 1 |

Correlation is significant at the 0.10* or 0.01*** level.

### 6. Conclusion and Future Research

In this study, we develop a methodological framework based on dynamic DEA that simultaneously estimates efficiency scores across multiple organizational dimensions while also generating an aggregate efficiency measure. The key findings and contributions of this research are summarized below:

- *Dual Modeling Approach:* We propose two complementary models: the first (SBM) estimates an overall efficiency score and then derives dimension-specific efficiencies, while the second (GP-SBM) first estimates dimension-level efficiencies and then computes an aggregate score.
- *Methodological Robustness:* Both models directly incorporate desirable and undesirable outputs without requiring data transformations, and employ a regularization parameter to enhance discrimination in multidimensional DEA settings and cases with many input and output variables.
- *Consistency and Discrimination:* Numerical analyses demonstrate that the regularization improves discrimination among DMUs and yields efficiency estimates that vary smoothly as its strength increases, i.e., there are no abrupt fluctuations or ranking reversals across hospitals.
- *Empirical Validation:* A case study of 12 Ontario hospitals evaluated over 24 months across three dimensions (technical, clinical, and patient experience), with shared input and output variables, confirms that both models capture dynamic, multidimensional efficiency patterns. Because our approach explicitly accounts for shared inputs/outputs and overlapping processes, the models provide a less biased assessment of efficiency rankings than conventional DEA-based approaches.

The proposed framework has several limitations. First, as a non-parametric frontier method, DEA attributes all deviations from the efficiency frontier to operational inefficiencies and does not explicitly account for statistical noise or measurement error. This limitation is particularly relevant in multidimensional settings, where there are numerous inputs, outputs, and carry-over variables that may introduce measurement variability. In such contexts, attributing all deviations from the frontier to inefficiency may be restrictive. Incorporating stochastic frontier approaches or bootstrapping techniques could provide additional insight into the statistical reliability of the estimated efficiency scores. A second limitation concerns the specification of efficiency dimensions and the assignment of variables to those dimensions. Although the proposed framework simultaneously evaluates multiple efficiency dimensions, the selection of which inputs and outputs belong to each must still be defined ex ante by a decision maker. As a result, alternative conceptualizations of organizational performance may lead to different model specifications and, consequently, different institutional rankings. Third, while the regularization mechanism improves discriminatory power when the number of variables is large relative to the number of DMUs, its value introduces an additional modeling decision. Although the sensitivity analysis indicates that rankings remain stable across a range of values, the parameter still represents an externally specified tuning choice rather than an endogenously determined quantity. Finally, GP-SBM relies on a linear approximation of an underlying nonlinear multi-objective formulation. Although the numerical experiments demonstrate strong correspondence between the approximation and the nonlinear model, the resulting efficiency scores should be interpreted as approximations to the original multi-objective problem rather than exact solutions.

Future research could extend this work in several directions. First, alternative regularization schemes or variable selection methods directly embedded within the dynamic DEA framework could further enhance discriminatory power in large-scale applications. Second, other multi-objective optimization approaches – such as ε-constraint or parametric methods – could be applied in place of goal programming to integrate multiple efficiency dimensions. Finally, the framework could be generalized to non-radial or range-adjusted formulations, enabling broader application in healthcare and other service systems where interdependent processes and shared resources create overlapping efficiency dimensions.

## References

Adler, N., & Golany, B. (2001). Evaluation of deregulated airline networks using data envelopment analysis combined with principal component analysis with an application to Western Europe. *European Journal of Operational Research*, 132(2), 260-273.

Afonso, G. P., Ferreira, D. C., & Figueira, J. R. (2023). A Network-DEA model to evaluate the impact of quality and access on hospital performance. *Annals of Operations Research*, 1-31.

Ali, A. I., & Seiford, L. M. (1990). Translation invariance in data envelopment analysis. O*perations Research Letters*, 9(6), 403-405.

Behzadian, M., Otaghsara, S. K., Yazdani, M., & Ignatius, J. (2012). A state-of the-art survey of TOPSIS applications. *Expert Systems with Applications*, *39*(17), 13051-13069.

Benıtez-Pen˜a, S., Bogetoft, P., and Morales, D. R. (2020). Feature selection in data envelopment analysis: A mathematical optimization approach. *Omega*, 96:102068.

Cantor, V. J. M., & Poh, K. L. (2018). Integrated analysis of healthcare efficiency: a systematic review. *Journal of Medical Systems*, 42, 1-23.

Chang, C. T. (2007). Multi-choice goal programming. *Omega*, *35*(4), 389-396.

Charnes, A., & Cooper, W. W. (1962). Programming with linear fractional functionals. *Naval Research logistics quarterly*, *9*(3-4), 181-186.

Charnes, A., Cooper, W. W., & Ferguson, R. O. (1955). Optimal estimation of executive compensation by linear programming. *Management Science*, *1*(2), 138-151.

Charnes, A., Cooper, W. W., & Rhodes, E. (1978). Measuring the efficiency of decision making units. *European journal of operational research*, *2*(6), 429-444.

Chen, S., Wu, Y., Chen, Y., Zhu, H., Wang, Z., Feng, D., ... & Liu, Z. (2016). Analysis of operation performance of general hospitals in Shenzhen, China: a super-efficiency data envelopment analysis. *The Lancet*, 388, S57.

Chen, Z. and Han, S. (2021). Comparison of dimension reduction methods for DEA under big data via Monte Carlo simulation. *Journal of Management Science and Engineering*, 6(4):363–376.

Cheng, G., & Zervopoulos, P. D. (2014). Estimating the technical efficiency of health care systems: A cross-country comparison using the directional distance function. *European Journal of Operational Research*, 238(3), 899-910.

Chung, Y. H., Färe, R., & Grosskopf, S. (1997). Productivity and undesirable outputs: a directional distance function approach. *Journal of Environmental Management*, 51(3), 229-240.

Dyson, R. G., Allen, R., Camanho, A. S., Podinovski, V. V., Sarrico, C. S., and Shale, E. A. (2001). Pitfalls and protocols in DEA. *European Journal of Operational Research*, 132(2):245–259.

Emrouznejad, A., Brzezicki, Ł., & Lu, C. (2025). Development and Evolution of Slacks-Based Measure Models in Data Envelopment Analysis: A Comprehensive Review of the Literature. *Journal of Economic Surveys*. https://doi.org/10.1111/joes.12682

Färe, R., Grosskopf, S., Lovell, C. K., & Pasurka, C. (1989). Multilateral productivity comparisons when some outputs are undesirable: a nonparametric approach. *The review of economics and statistics*, 90-98.

Fazria, N. F., & Dhamayanti, I. (2021). A literature review on the identification of variables for measuring hospital efficiency in the data envelopment analysis (DEA). *Telemedicine Use in Health Facility During Covid-19 Pandemic: Literature Review*, 10(1), 1-15.

Ferreira, D. C., Figueira, J. R., Greco, S., & Marques, R. C. (2023). Data envelopment analysis models with imperfect knowledge of input and output values: An application to Portuguese public hospitals. *Expert Systems with Applications*, 231, 120543.

Fukuyama, H., Tsionas, M., & Tan, Y. (2023). Dynamic network data envelopment analysis with a sequential structure and behavioral-causal analysis: application to the Chinese banking industry. *European Journal of Operational Research*, *307*(3), 1360-1373.

Gong, J., Shi, L., Wang, X., & Sun, G. (2023). The efficiency of health resource allocation and its influencing factors: evidence from the super efficiency slack based model-Tobit model. *International health*, 15(3), 326-334.

Gupta, S., Bandyopadhyay, G., Bhattacharjee, M., & Biswas, S. (2019). Portfolio Selection using DEA-COPRAS at risk–return interface based on NSE (India). *International Journal of Innovative Technology and Exploring Engineering (IJITEE)*, *8*(10), 4078-4086.

Halkos, G., & Petrou, K. N. (2019). Treating undesirable outputs in DEA: A critical review. *Economic Analysis and Policy*, 62, 97-104.

Hua, Z., Bian, Y. (2007). DEA with Undesirable Factors. In: Zhu, J., Cook, W.D. (eds) Modeling Data Irregularities and Structural Complexities in Data Envelopment Analysis. Springer, Boston, MA.

Kao, C., Pang, R. Z., Liu, S. T., & Bai, X. J. (2021). Most productive types of hospitals: An empirical analysis. *Omega (United Kingdom)*, 99, 102310.

Lee, C.-Y., & Cai, J.-Y. (2020). Lasso variable selection in data envelopment analysis with small datasets. *Omega*, 91, 102019.

Lin, F., Deng, Y. J., Lu, W. M., & Kweh, Q. L. (2019). Impulse response function analysis of the impacts of hospital accreditations on hospital efficiency. *Health Care Management Science*, 22, 394-409.

Liu, W., & Sharp, J. (1999). DEA models via goal programming. *Data envelopment analysis in the service sector*, 79-101.

Lovell, C. K., Pastor, J. T., & Turner, J. A. (1995). Measuring macroeconomic performance in the OECD: A comparison of European and non-European countries. *European Journal of Operational Research*, 87(3), 507-518.

Mahdiloo, M., Saen, R. F., & Lee, K. H. (2015). Technical, environmental and eco-efficiency measurement for supplier selection: An extension and application of data envelopment analysis. *International journal of production economics*, *168*, 279-289.

Mavrotas, G. (2009). Effective implementation of the ε-constraint method in multi-objective mathematical programming problems. *Applied Mathematics and Computation*, *213*(2), 455-465.

Mergoni, A., Emrouznejad, A., & De Witte, K. (2025). Fifty years of data envelopment analysis. *European Journal of Operational Research*, 326(3), 389-412.

Nataraja, N. R. and Johnson, A. L. (2011). Guidelines for using variable selection techniques in data envelopment analysis. *European Journal of Operational Research*, 215(3):662–669.

Omrani, H. (2013). Common weights data envelopment analysis with uncertain data: A robust optimization approach. *Computers & Industrial Engineering*, 66(4), 1163-1170.

Omrani, H., Shamsi, M., & Emrouznejad, A. (2022). Evaluating sustainable efficiency of decision-making units considering undesirable outputs: an application to airline using integrated multi-objective DEA-TOPSIS. *Environment, Development and Sustainability*, *25*(7), 5899-5930.

Ortega-Díaz, M. I., & Martín, J. C. (2022). How to detect hospitals where quality would not be jeopardized by health cost savings? A methodological approach using DEA with SBM analysis. *Health Policy*, 126(10), 1069-1074.

Pai, D. R., Pakdil, F., & Azadeh-Fard, N. (2024). Applications of data envelopment analysis in acute care hospitals: a systematic literature review, 1984–2022. Health Care Management Science, 1-29.

Paradi, J. C., Rouatt, S., & Zhu, H. (2011). Two-stage evaluation of bank branch efficiency using data envelopment analysis. *Omega*, *39*(1), 99-109.

Pastor, J. T., Ruiz, J. L., and Sirvent, I. (2002). A statistical test for nested radial DEA models. *Operations Research*, 50(4):728–735.

Peixoto, M. G. M., Musetti, M. A., & de Mendonça, M. C. A. (2020). Performance management in hospital organizations from the perspective of Principal Component Analysis and Data Envelopment Analysis: the case of Federal University Hospitals in Brazil. *Computers & Industrial Engineering*, 150, 106873.

Qin, Z. and Song, I. (2014). Joint variable selection for data envelopment analysis via group sparsity. *arXiv preprint arXiv*:1402.3740.

Rakhshan, S. A. (2017). Efficiency ranking of decision making units in data envelopment analysis by using TOPSIS-DEA method. *Journal of the Operational Research Society*, *68*(8), 906-918.

Roth, A., Tucker, A. L., Venkataraman, S., & Chilingerian, J. (2019). Being on the productivity frontier: Identifying "triple aim performance" hospitals. *Production and Operations Management*, 28(9), 2165-2183.

Ruggiero, J. (2005). Impact assessment of input omission on DEA. *International Journal of Information Technology & Decision Making*, 4(03):359–368.

Saaty, T., & Gass, S. (1954). Parametric objective function (part 1). *Journal of the Operations Research Society of America*, *2*(3), 316-319.

Simar, L. and Wilson, P. W. (2001). Testing restrictions in nonparametric efficiency models. *Communications in Statistics-Simulation and Computation*, 30(1):159–184.

Sommersguter-Reichmann, M. (2022). Health care quality in nonparametric efficiency studies: a review. *Central European Journal of Operations Research*, 30(1), 67-131.

Sun, M., Ye, Y., Zhang, G., Shang, X., & Xue, Y. (2023). Healthcare services efficiency and its intrinsic drivers in China: based on the three-stage super-efficiency SBM model. *BMC Health Services Research*, 23(1), 811.

Tone, K. (2001). A slacks-based measure of efficiency in data envelopment analysis. *European Journal of Operational Research*, 130(3), 498-509.

Tone, K., & Tsutsui, M. (2010). Dynamic DEA: A slacks-based measure approach. *Omega*, *38*(3-4), 145-156.

Verma, U., Diamant, A., Imanirad, R., Verma, A., & Razak, F. (2025). Resilience and Adaptation: Assessing the Effects of COVID-19 on Hospital Efficiency and Quality. *Available at SSRN 5396361*.

Wagner, J. M. and Shimshak, D. G. (2007). Stepwise selection of variables in data envelopment analysis: Procedures and managerial perspectives. *European Journal of Operational Research*, 180(1):57–67.

Wang, M., Chen, Y., & Zelenyuk, V. (2025). DEA for big wide data based on regularization approaches: An application to energy efficiency analysis. *Energy Economics*, 108919.

Wang, Z., Yang, M., Liang, L., & Feng, C. (2026). Integrating LASSO and DEA via Nash bargaining–inspired framework for variable selection under high dimensional data. *Omega*, 103551.

Yang, F., Wei, F., Li, Y., Huang, Y., & Chen, Y. (2018). Expected efficiency based on directional distance function in data envelopment analysis. *Computers & Industrial Engineering*, 125, 33-45.

Yilmaz, I. (2023). A hybrid DEA–fuzzy COPRAS approach to the evaluation of renewable energy: A case of wind farms in Turkey. *Sustainability*, *15*(14), 11267.

Zarrin, M. (2023). A mixed-integer slacks-based measure data envelopment analysis for efficiency measuring of German university hospitals. *Health Care Management Science*, 26(1), 138-160.

Zavadskas, E. K., Kaklauskas, A., Peldschus, F., & Turskis, Z. (2007). Multi-attribute assessment of road design solutions by using the COPRAS method. *The Baltic of Road and Bridge Engineering*, *2*(4), 195-203.

Zhang, X., Tone, K., & Lu, Y. (2018). Impact of the local public hospital reform on the efficiency of medium-sized hospitals in Japan: An improved slacks-based measure data envelopment analysis approach. *Health Services Research*, 53(2), 896-918.